%% file: arxiv2.tex
\documentclass[smallextended,numbook,runningheads]{svjour3} 
\smartqed
\usepackage[whole]{bxcjkjatype}

\usepackage{lscape}
\usepackage{pdflscape}
\usepackage{lipsum} 
\usepackage{geometry}
\usepackage{graphicx}
\usepackage{amsmath}
\usepackage{mathptmx}
\usepackage{amsfonts}
\usepackage{color}
\usepackage{xspace}
\usepackage{array}
\usepackage{pgfplots}
\usepackage{ifthen}
\usepackage{mathtools}
\mathtoolsset{centercolon}
\usepackage{algorithm}
\usepackage{algorithmic}
\usepackage[colorlinks=true,pagebackref=true,citecolor = cyan,linkcolor=cyan,hypertexnames=false]{hyperref}
\renewcommand*{\backref}[1]{}
\renewcommand*{\backrefalt}[4]{[{\footnotesize%
		\ifcase #1 Not cited.%
		\or Cited on page~#2.%
		\else Cited on pages #2.%
		\fi%
	}]}
\usepackage{cleveref}
\usepackage{autonum}
\pgfplotscreateplotcyclelist{mylist}{%
teal,every mark/.append style={fill=teal!80!black},mark=star\\%
orange,every mark/.append style={fill=orange!80!black},mark=*\\%
cyan!60!black,every mark/.append style={fill=cyan!80!black},mark=otimes\\%
red,every mark/.append style={solid,fill=red!80!black},mark=pentagon\\%
magenta,mark=square\\%
lime!60!black,every mark/.append style={fill=lime},mark=diamond\\%
yellow!60!black,densely dashed,
every mark/.append style={solid,fill=yellow!80!black},mark=star\\%
black,every mark/.append style={solid,fill=gray},mark=*\\%
blue,densely dashed,mark=otimes,every mark/.append style=solid\\%
red,densely dashed,every mark/.append style={solid,fill=red!80!black},mark=pentagon\\%
}
\usepackage{tikz}
\usepackage{pgfplotstable}
\usetikzlibrary{arrows,positioning,plotmarks,external,patterns,angles,
decorations.pathmorphing,backgrounds,fit,shapes,calc}

\DeclarePairedDelimiterX{\setI}[2]{\{}{\}}{\,{#1}\ \delimsize| \ {#2}\,}

\input{macro.tex}
\newcommand{\hA}{\hat{A}}
\newcommand{\hB}{\hat{B}}
\newcommand{\hC}{\hat{C}}
\newcommand{\hD}{\hat{D}}

\crefname{theorem}{Theorem}{Theorems}
\crefname{lemma}{Lemma}{Lemmas}
\crefname{proposition}{Proposition}{Propositions}
\crefname{corollary}{Corollary}{Corollaries}
\crefname{definition}{Definition}{Definitions}
\crefname{example}{Example}{Examples}
\crefname{remark}{Remark}{Remarks}
\crefname{figure}{Figure}{Figures}
\crefname{table}{Table}{Tables}

\crefname{section}{Section}{Sections}
\crefname{subsection}{Section}{Sections}

\newcounter{NN}
\newcounter{NR}
\newcounter{ncount}
\newcounter{scount}

\def\black{\setcounter{NN}{1}}
\def\white{\setcounter{NN}{2}}
\def\change{\ifthenelse{\value{NN}=1}%
{\setcounter{NN}{2}}%
{\setcounter{NN}{1}}}

\def\metrics#1#2#3#4{%
\def\treesize{#1}\def\thick{#2}\def\rad{#3}%
\pgfmathparse{#3-#4}\let\smallrad=\pgfmathresult}

\def\vertex#1#2{%
\addtocounter{ncount}{1}
\setcounter{scount}{\thencount}
\addtocounter{scount}{-#1}
\ifthenelse{\first=1}{%
\node (\thencount) at ($(\thescount) + ( #2,1)$){};
\draw (\thencount)--(\thescount);
\fill (\thencount) circle;
}%
{
\ifthenelse{\value{NN}=2}{\fill[white,radius=\smallrad pt] (\thencount) circle;}{}
}}

\def\tree#1{%
\ifthenelse{\value{NN}=0}{}{\setcounter{NR}{\value{NN}}}
\begin{tikzpicture}[x=\treesize mm,y=\treesize mm,radius=\rad pt,line width=\thick pt,inner sep=0,baseline=-0.02cm]
\node (0) at (0,0) {}; \fill (0) circle;
\setcounter{NN}{\value{NR}}
\setcounter{ncount}{0}
\gdef\first{1}
#1
\setcounter{NN}{\value{NR}} 
\gdef\first{2}
\ifthenelse{\value{NR}=2}{\fill[white,radius=\smallrad pt] (0) circle;}{}
\setcounter{ncount}{0}
#1
\setcounter{NN}{0}
\end{tikzpicture}}

\journalname{}

\begin{document}

\title{A new family of fourth-order energy-preserving integrators}


\author{Yuto Miyatake}


\institute{
        Yuto Miyatake \at
        Cybermedia Center,
        Osaka University, Osaka, Japan.
        \email{yuto.miyatake.cmc@osaka-u.ac.jp}
}

\date{Received: date / Accepted: date}

\maketitle

\metrics{2.5}{1}{2}{0.8}

\begin{abstract}{
For Hamiltonian systems with non-canonical structure matrices,  a new family of fourth-order energy-preserving integrators is presented. 
The integrators take a form of a combination of Runge--Kutta methods and continuous-stage Runge--Kutta methods and feature a set of free parameters that offer greater flexibility and efficiency. Specifically, we demonstrate that by carefully choosing these free parameters, a simplified Newton iteration applied to the integrators of order four can be parallelizable. 
This results in faster and more efficient integrators compared with existing fourth-order energy-preserving integrators.
}
\end{abstract}

\section{Introduction}
\label{sec:intro}
In this paper, we are concerned with the numerical integration of a system of ordinary differential equations (ODEs) of the form
\begin{equation}
    \label{eq:poisson}
    \frac{\rmd}{\rmd t} y = S(y) \nabla H(y), \quad y(0) = y_0 \in \bbR^d,
\end{equation}
where $y:[0,T)\to\bbR^d$ is a dependent variable, $S:\bbR^d \to \bbR^{d\times d}$ is a skew-symmetric matrix function, and $H:\bbR^d \to \bbR$ is a real-valued function, which we call energy. 
The two functions $S$ and $H$ are assumed to be sufficiently regular.
Along the solution to \eqref{eq:poisson}, the function $H$ is constant:
\begin{equation}
    \frac{\rmd}{\rmd t} H(y(t)) = \nabla H(y(t)) ^\top \dot{y}(t) = \nabla H(y(t)) ^\top S(y(t))  \nabla H(y(t)) = 0,
\end{equation}
where the dot stands for the differentiation with respect to $t$.
Conversely, a system of ODEs having a first-integral can always be formulated in the form of \eqref{eq:poisson} with an appropriate skew-symmetric matrix function $S(y)$~\cite{qc96},
although the expression of $S(y)$ might not be unique.
When $S(y)$ is constant and in particular $S = J^{-1}$ with
\begin{equation}
    J = \begin{bmatrix}
    O & -I \\ I & O
    \end{bmatrix},
\end{equation}
the system is called a Hamiltonian system and the corresponding $H$ is often referred to as the Hamiltonian. 
In more general cases, where $S$ depends on $y$, if the Poisson bracket satisfies the Jacobi identity, the system is called a Poisson system (see, e.g., \cite[Chapter~VII.2]{hlw06}).
In this paper, we always call the system of the form \eqref{eq:poisson} a Poisson system, even if the Poisson bracket does not satisfy the Jacobi identity, bearing in mind that this terminology is just for convenience only.
Furthermore, depending on the structure of $S(y)$, the system of the form \eqref{eq:poisson} exhibits rich geometric properties, such as symplecticity.

This paper focuses on energy-preserving integrators,  which are a typical branch of geometric numerical integrators~\cite{hlw06}.
In this paper, an energy-preserving integrator refers to a one-step method $y_0\mapsto y_1$, where $y_1\approx y(h)$, such that $H(y_1)=H(y_0)$.
For such a method, $H(y_n)=H(y_0)$ holds even if the step-size $h$ is controlled adaptively.

The projection method is a relatively simple method.
The projection method, while conceptually straightforward, encompasses a variety of approaches for projecting solutions onto the appropriate manifold.
The effectiveness of this method, particularly concerning long-term behaviour, may depend on both the selected projection technique and the underlying integrator~\cite[Chapter~IV.4]{hlw06}.
Therefore, caution is advised in employing the projection method.
\footnote{It is worth noting that the projection concept retains its utility. Additionally, we mention the work~\cite{nm15} that shows the equivalence between projection methods and the discrete gradient methods, which will mentioned below.}
A more sophisticated and systematic approach is called the discrete gradient method, which was first formulated by Gonzalez~\cite{go96} (see also McLachlan, Quispel and Robidoux~\cite{mq99}), although a similar idea had been known for quite some time.
The discrete gradient method usually produces a second-order energy-preserving integrator.
The average vector field (AVF) method, proposed by Quispel and McLaren~\cite{qm08}, is a subclass of the discrete gradient method, which is a B-series method and of order two when $S$ is a constant skew-symmetric matrix. 
Over the past decade, extensions of the AVF method to higher order have been extensively studied.
For a constant $S$, Hairer proposed the AVF collocation method~\cite{ha10} and Brugnano, Iavernaro and Trigiante proposed the Hamiltonian boundary value method~\cite{br10}.
These methods are based on so-called continuous-stage Runge--Kutta methods.
Also worth mentioning is a relatively new work~\cite{ei22}, which establishes a general theory on the order theory for discrete gradient methods.

Roughly speaking, the computational cost of the AVF collocation method is almost the same as the Gauss method of the same order.
Miyatake and Butcher~\cite{mb16} characterize the condition for the continuous stage Runge--Kutta (CSRK) method being energy-preserving in terms of the symmetry of an $s\times s$ matrix $M$ defining the CSRK method and find that the order condition can also be characterized in terms of $M^{-1}$.
The discussions using $M$ seem fruitful in that 
one can construct an integrator with an intended order with some degrees of freedom.
By
manipulating the remaining parameters to enhance the integrator, for example, parallelizable integrators can be constructed.

For Poisson systems, 
efforts developing energy-preserving methods have also been devoted, and the aforementioned methods have been extended to this general class.
For example, the AVF collocation method is extended to Poisson systems by introducing a new class of integration methods, which is a generalization of the CSRK method to partitioned systems~\cite{ch11}
(see~\cite{am22,br12} for the extension of the Hamiltonian boundary value method to the Poisson systems).
We refer to this new class of integration method as the partitioned continuous-stage Runge--Kutta (PCSRK) method.
The $s$-degree PCSRK method\footnote{The degree plays a similar role to the stages of (partitioned) Runge--Kutta methods.} is characterized by the $s\times s$ matrices $M_i$ ($i=1,\dots,s$) and the nodes $c_i$ ($i=1,\dots,s$).
The results given in~\cite{mi15} suggest that a PCSRK method is energy-preserving if all $M_i$ are symmetric, but it is not clear if the order condition is concisely characterized in terms of the matrices $M_i$ and nodes $c_i$, as is clearly done for the CSRK methods.
This task does not seem so trivial; there are several difficulties associated with it. 
For example, recall that for the constant $S$, the order condition is characterized in terms of the \emph{inverse} of $M$; however, for Poisson systems, although the highest order of the $s$-degree PCSRK methods is $2s$~\cite{ch11}, the corresponding $M_i$'s are singular. 
Thus, one cannot expect that the order conditions are characterized in terms of the inverse of $M_i$'s.
Other difficulties are discussed in \cref{rem:difficulty}.

Taking the above backgrounds into consideration, we focus only on fourth-order methods and address the following issues.
\begin{itemize}
    \item The fourth-order PCSRK method exists with $s=2$~\cite{ch11}, which is unique if the degree is restricted to $s=2$.
    In this paper, we set $s=3$ and characterize the method for being order 4 in terms of $M_1,M_2,M_3$ and $c_1,c_2,c_3$.
    The key idea is to require that the PCSRK method be reduced to the CSRK method discussed in~\cite{mb16} when $S$ is constant, and the method is symmetric, and simplify the order conditions for the bi-coloured rooted trees with three vertices having a black root. 
    This can be viewed as a three-degree PCSRK method with some degrees of freedom, and the parameters can be devised from another perspective.
    \item As discussed in~\cite{mb16}, an advantage of a numerical method with some degrees of freedom is that much more efficient variants may be able to be explored. 
    Clearly, in general, the larger the degrees of the CSRK methods are, the more expensive the computational cost becomes. 
    However, this is not always the case; for example, if the matrix $M$ has a specific eigenstructure, the method can be implemented in a parallel architecture with almost the same cost as the case $s=1$, though the memory usage grows.
    We show that a similar structure holds for the PCSRK methods.
    \item Based on the above two points, we develop three-degree PCSRK integrators with some degrees of freedom, which are energy-preserving for Poisson systems, of order four, and efficiently implemented in a parallel architecture.
    The proposed integrators are reduced to the ones developed in~\cite{mb16} with similar properties for Hamiltonian systems.
\end{itemize}

The paper is organized as follows.
In Section~\ref{sec:2}, after reviewing energy-preserving CSRK methods for constant $S$ and their properties, we also discuss the formulation of energy-preserving PCSRK methods for general cases.
We develop a family of energy-preserving PCSRK integrators in Section~\ref{sec:3}. We discuss the implementation issues and optimal parameter choices in Section~\ref{sec:4}.
Concluding remarks are given in Section~\ref{sec:5}.

\section{Preliminaries}
\label{sec:2}

\subsection{Hamiltonian systems and CSRK methods}

Let us consider the system
\begin{align}
    \label{eq:hamilton}
    \frac{\rmd}{\rmd t} y = S \nabla H(y), \quad y(0) = y_0 \in \bbR^d,
\end{align}
where $S$ is a constant skew-symmetric matrix, but not necessarily $J^{-1}$. 
Hamiltonian systems are a typical example of this class.
The average vector field (AVF) method reads
\begin{align}
    \label{eq:avf}
    y_1 = y_0 + h\int_0^1 f((1-\tau) y_0 + \tau y_1) \,\rmd \tau,  
\end{align}
where $f(y) = S \nabla H(y)$.
This method is of second order and energy-preserving $H(y_1)=H(y_0)$.
The method can be regarded as a continuous stage Runge--Kutta method.

\begin{definition}[CSRK methods]
Let $A_{\tau,\zeta}$ be a polynomial in $\tau$ and $\zeta$.
Assume that $A_{0,\zeta}=0$.
The polynomial degree of $A_{\tau,\zeta}$ in $\tau$ is denoted by $s$.
Let $B_\zeta$ be defined by $B_\zeta = A_{1,\zeta}$.
Define an $s$-degree polynomial $Y_\tau$ ($\tau\in[0,1]$)
and $y_1$ such that they satisfy
\begin{align}
    Y_\tau &= y_0 + h \int_0^1 A_{\tau,\zeta} f(Y_\zeta)\,\rmd \zeta,\\
    y_1 &= y_0 + h \int_0^1 B_\tau f(Y_\tau)\,\rmd \tau.
\end{align}
A one-step method $y_0\mapsto y_1$ is called an $s$-degree continuous stage Runge--Kutta (CSRK) method.
\end{definition}

The above definition does not specify the order of $A_{\tau,\zeta}$ in terms of $\zeta$, but
let us focus on $A_{\tau,\zeta}$ which is a polynomial of degree $s$ in $\tau$ and $s-1$ in $\zeta$.
Such $A_{\tau,\zeta}$ will be denoted by
\begin{equation}
    A_{\tau,\zeta}
    =
    \begin{bmatrix}
    \tau & \frac{\tau^2}{2} & \cdots & \frac{\tau^s}{s}
    \end{bmatrix}
    M
    \begin{bmatrix}
    1 \\ \zeta \\ \vdots \\ \zeta^{s-1}
    \end{bmatrix}
\end{equation}
with a constant matrix $M\in\bbR^{s\times s}$
so that the polynomial is identified with the matrix $M$.
When $s=1$ and $M=1$, the method reduces to the AVF method.

A sufficient condition for energy-preservation can be characterized in terms of $M$.

\begin{theorem}[\cite{mi14,mb16}, see also~\cite{ta14}]
When applied to \eqref{eq:hamilton},
a CSRK method is energy-preserving if $M$ is symmetric.
\end{theorem}

The symmetry of $M$ means $(\partial/\partial \tau)A_{\tau,\zeta}$ is symmetric.
The condition is also necessary under a mild condition~\cite{mb16}.

Several characterizations of the order conditions with respect to $M$ are given in~\cite{mb16}.
We note that the discussion there is based on the simplifying assumptions (see also~\cite{ha10}).
A CSRK method is energy-preserving and of order at least $p=2\eta$ if the symmetric matrix $M\in\bbR^{s\times s}$ satisfies
\begin{equation}
    \begin{bmatrix}
    \tfrac{1}{k} & \tfrac{1}{k+1} & \cdots & \tfrac{1}{k+s-1}
    \end{bmatrix}
    M
    = i_k^\trans, \quad k = 1,\dots,\eta,
\end{equation}
where $i_k$ denotes the $k$-th column of the $s\times s$ identity matrix.
If we choose $\eta = s$ the method with 
\begin{equation}
\renewcommand{\arraystretch}{1.2}
    \begin{bmatrix}
    1 & \tfrac{1}{2} & \cdots & \tfrac{1}{s} \\
    \tfrac{1}{2} & \tfrac{1}{3} & \cdots & \tfrac{1}{s+1} \\
    \vdots & \vdots & \ddots & \vdots \\
    \tfrac{1}{s} & \tfrac{1}{s+1} & \cdots & \tfrac{1}{2s-1}
    \end{bmatrix}
    M = I_s
\end{equation}
is of order $2s$ and coincides with the AVF collocation method of order $2s$~\cite{ha10}.
For example, when $s=2$, we have
\begin{equation}
    \label{eq:AVF_collocation_2}
    M = \begin{bmatrix}
    4 & -6 \\ -6 & 12
    \end{bmatrix}.
\end{equation}
As another illustrative example,
when $s=3$, the above characterization indicates that
the method with
\begin{equation}
\label{eq:M3alpha}
{\renewcommand{\arraystretch}{1.2}
\begin{bmatrix}
1 & \tfrac{1}{2} & \tfrac{1}{3} \\
\tfrac12 & \tfrac13 & \tfrac14 \\
\tfrac13 & \tfrac14 & \alpha
\end{bmatrix}
M = I_3, \quad \alpha \in \bbR
}
\end{equation}
is of order four as $\eta=2$ except for $\alpha = 1/5$.
In this way, one can construct high-order energy-preserving integrators with some degrees of freedom.
When $\alpha \neq 7 / 36$, by introducing a new variable (parameter) $\tilde{\alpha} = 1/(36\alpha - 7)$ the matrix $M$ can be expressed more explicitly as
{\renewcommand{\arraystretch}{1.2}
\begin{equation}
    M = \begin{bmatrix}
    \tilde{\alpha} + 4 & -6\tilde{\alpha} - 6 & 6 \tilde{\alpha} \\
    -6\tilde{\alpha} - 6 & 36 \tilde{\alpha} + 12 & -36 \tilde{\alpha} \\
    6\tilde{\alpha} & -36 \tilde{\alpha} & 36 \tilde{\alpha}
    \end{bmatrix}.
\end{equation}
}

\subsection{Poisson systems and PCSRK methods}

Note that the AVF method~\eqref{eq:avf} is not energy-preserving when applied to~\eqref{eq:poisson} with a non-constant $S(y)$.
A straightforward modification
\begin{equation} \label{eq:avf_extention}
    y_1 = y_0 + S \Big( \frac{y_1+y_0}{2} \Big) \int_0^1 \nabla H((1-\tau) y_0 + \tau y_1) \,\rmd \tau
\end{equation}
is energy-preserving, symmetric and thus of order two.
Here, $S(y)$ is discretized by using the mid-point rule to ensure the method is symmetric.
Other choices, such as $S(y_0)$ and $S(y_1)$, still guarantee the energy-preservation, though the resulting integrator is of order one.
It should be noted that in \eqref{eq:avf_extention} $\nabla H(y)$ term is discretized in a CSRK manner while $S(y)$ in a standard RK manner.
This observation leads to the following class of numerical integrators applied to a partitioned system
\begin{align}
    \frac{\rmd}{\rmd t} y &= S(z) \nabla H(y), \quad y(0) = y_0 , \\
    \frac{\rmd}{\rmd t} z &= S(z) \nabla H(y), \quad z(0) = z_0.
\end{align}

\begin{definition}[PCSRK methods]
\label{def:PCSRK}
Let $A_{i,\tau,j,\zeta}$ $(j=1,\dots,s)$ be a polynomial in $\tau$ and $\zeta$ with the property $A_{i,0,j,\zeta}=0$.
$A_{i,\tau,j,\zeta}$ is assumed to be independent of $i$; thus, it is often denoted by $A_{\tau,j,\zeta}$.
$\hA_{i,\tau,j,\zeta}$ is defined by $A_{i,c_i,j,\zeta}$ with $s$ distinct nodes $0\leq c_1< \dots< c_s\leq 1$.
Let $B_{j,\zeta} = \hB_{j,\zeta} = A_{1,j,\zeta}$.
Define an $s$-degree polynomial $Y_\tau$ ($\tau\in[0,1]$),
$Z_1,\dots,Z_s$, $y_1$ and $z_1$ such that they satisfy
\begin{align}
    Y_\tau &= y_0 + h \sum_{j=1}^s \int_0^1 A_{i,\tau,j,\zeta} S(Z_j) \nabla H (Y_\zeta) \,\rmd \zeta,\label{pcsrk1}\\
    Z_i &= z_0 + h \sum_{j=1}^s \int_0^1 \hA_{i,\tau,j,\zeta} S(Z_j) \nabla H (Y_\zeta) \,\rmd \zeta, \label{pcsrk2}\\
    y_1 &= y_0 + h \sum_{i=1}^s \int_0^1 B_{i,\tau} S(Z_i) \nabla H (Y_\tau) \,\rmd \tau,\label{pcsrk3}\\
    z_1 &= z_0 + h \sum_{i=1}^s \int_0^1 \hB_{i,\tau} S(Z_i) \nabla H (Y_\tau) \,\rmd \tau \label{pcsrk4}
\end{align}
with $y_0=z_0$.
A one-step method $y_0\mapsto y_1$ is called an $s$-degree partitioned CSRK (PCSRK) method.
\end{definition}

We note that by definition $Z_i = Y_{c_i}$ and $y_1 = z_1$; thus, the scheme can be written in a more compact form
\begin{align}
    Y_\tau &= y_0 + h \sum_{j=1}^s \int_0^1 A_{i,\tau,j,\zeta} S(Y_{c_j}) \nabla H (Y_\zeta) \,\rmd \zeta,\label{pcsrk1_comp} \\
    y_1 &= y_0 + h \sum_{i=1}^s \int_0^1 B_{i,\tau} S(Y_{c_i}) \nabla H (Y_\tau) \,\rmd \tau.\label{pcsrk3_comp}
\end{align}
The expression in \cref{def:PCSRK} is useful for discussing the order conditions.

As $A_{i,\tau,j,\zeta}$ depends on $\tau$, $j$ and $\zeta$, but does not depend on $i$, it is convenient to express it as
\begin{equation}
    A_{i,\tau,j,\sigma}
    =
    \begin{bmatrix}
    \tau & \frac{\tau^2}{2} & \cdots & \frac{\tau^s}{s}
    \end{bmatrix}
    M_j
    \begin{bmatrix}
    1 \\ \zeta \\ \vdots \\ \zeta^{s-1}
    \end{bmatrix}
\end{equation}
by using constant matrices $M_j\in\bbR^{s\times s}$ for $j=1,\dots,s$.
For example, the second-order method proposed by Cohen and Hairer~\cite{ch11} is given by
\begin{equation}
    M_1 =
    \begin{bmatrix}
    2+\sqrt{3} & -(3+\sqrt{3}) \\
    - (3+\sqrt{3}) & 6
    \end{bmatrix},
    \quad 
    M_2 =
    \begin{bmatrix}
    2-\sqrt{3} & \sqrt{3}-3 \\
    \sqrt{3}-3 & 6
    \end{bmatrix}
\end{equation}
and $c_1,c_2 = 1/2 \mp \sqrt{3}/6$.
We observe that
\begin{equation}
    M_1 + M_2 = \begin{bmatrix}
    4 & -6 \\ -6 & 12
    \end{bmatrix},
\end{equation}
which coincides with $M$ in \eqref{eq:AVF_collocation_2}.
Thus, the method is reduced to the AVF collocation method~\cite{ha10} when $S$ is constant.
It should be noted that both $M_1$ and $M_2$ are singular while $M_1+M_2$ is nonsingular.
This indicates that one cannot expect that $M_i$ is characterized as the inverse of some matrices.

Sufficient conditions of PCSRK methods to be energy-preserving are characterized in terms of $A_{i,\tau,j,\zeta}$, or equivalently, the $M_i$ matrices.

\begin{theorem}[\cite{mi15}]
\label{thm:pcsrk_energy}
When applied to \eqref{eq:poisson},
a PCSRK method is energy-preserving if all $M_i$'s are symmetric
\end{theorem}

The conditions for symmetry are also characterized as follows.

\begin{theorem}[cf.~\cite{ch11}]
\label{thm:pcsrk_symmetry}

A PCSRK method is symmetric if $A_{i,1-\tau,s+1-j,1-\zeta} + A_{i,\tau,j,\zeta} = B_{j,\zeta}$ and $c_{s+1-i}=1-c_i$.
\end{theorem}

Order conditions will be discussed in the next section.

\section{A family of fourth-order energy-preserving integrators}
\label{sec:3}

In this section, we derive a family of fourth-order energy-preserving PCSRK methods.
Specifically, we set $s=3$ and aim to derive a 3-degree PCSRK integrator with some free parameters.

We begin with a few remarks on order conditions.
In the formulation of the PCSRK method (\cref{def:PCSRK}), it is not necessary to calculate $z_1$, because
only $y_1$ is required as an output to proceed with the subsequent steps. The expression for $y_1$ is given as a P-series:
\begin{equation}\label{eq:PRK:P-series}
    y_1 = y_0 + \sum_{ \tau \in TP_y} \frac{ h^{|\tau|} }{\sigma (\tau)} \phi (\tau) F(\tau) (y_0,y_0),
\end{equation}
where
$\phi, \sigma$ and $F$ are the elementary weights, symmetry and elementary differentials, respectively.
The symbol $ TP_y $ denotes the set of bi-coloured trees with black roots,
i.e., 
\begin{equation}
    TP_y = \{ \black\tree{},
    \black\tree{\vertex10},
    \black\tree{\change\vertex10},
    \black\tree{\vertex1{0.6}\vertex2{-0.6}},
    \black\tree{\white\vertex1{0.6}\black\vertex2{-0.6}},
    \black\tree{\white\vertex1{0.6}\white\vertex2{-0.6}},
    \tree{\vertex10\vertex10},
    \tree{\vertex10\change\vertex10},
    \tree{\change\vertex10\change\vertex10},
    \tree{\change\vertex10\vertex10},\dots\}.
\end{equation}
For more details about the order conditions with bi-coloured trees, refer to~\cite{hlw06}.

A PCSRK method is of order $p$ if it satisfies
\begin{equation}
    \phi (\tau) = e (\tau)
    \quad \text{for} \quad
    \tau \in TP_y, \quad |\tau| \leq p,
\end{equation}
where $|\tau|$ denotes the order of $\tau$, i.e. the number of vertices of $\tau$.
Here, $e(\tau)$ for the mono-coloured trees is defined by
\begin{equation}
    e(\emptyset) = e(\black\tree{}) = 1,
    \quad
    e(\tau) = 
    \frac{1}{|\tau|} e(\tau_1) \cdots e(\tau_m)
    \quad \text{for}
    \quad
    \tau = [\tau_1,\dots,\tau_m]
\end{equation}
and $e(\tau)$ for the bi-coloured trees takes the same value with the mono-coloured trees.

\begin{remark}
\label{rem:difficulty}
Since checking all the order conditions is an immense task, one approach to avoid checking every condition relies on the simplifying assumptions. 
For the PCSRK methods, the simplifying assumptions are given as follows:
\begin{alignat}{2}
    B(\rho) &\quad  \sum_{i=1}^s \int_0^1 B_{i,\tau}C_{i,\tau}^{k-1} \hC_i^l \,\rmd \tau = \frac{1}{k+l}, & & \quad  1\leq k + l \leq \rho,  \\
    C(\eta) & \quad  \sum_{j=1}^s \int_0^1 A_{i,\tau,j,\sigma} C_{j,\sigma}^{k-1}\hC_j^l \, \rmd \sigma = \frac{1}{k+l} C_{i,\tau}^{k+l}, & & \quad  1\leq k + l \leq \eta, \\
    \hC (\eta) &\quad \sum_{j=1}^s \int_0^1 \hA_{i,\tau,j,\sigma} C_{j,\sigma}^{k-1}\hC_j^l \, \rmd \sigma = \frac{1}{k+l} \hC_{i,\tau}^{k+l}, & &\quad 1\leq k + l \leq \eta,\\
    D(\xi) &\quad \sum_{i=1}^s \int_0^1 B_{i,\tau} C_{i,\tau}^{k-1} \hC_i^l A_{i,\tau,j,\sigma} \, \rmd \tau = \frac{B_{j,\sigma}}{k+l}(1-\hC_{j,\sigma}^{k+l}), & &\quad 1 \leq k + l \leq \xi, \\
    \hD(\xi) &\quad \sum_{i=1}^s \int_0^1 B_{i,\tau} C_{i,\tau}^{k-1} \hC_i^l \hA_{i,\tau,j,\sigma} \, \rmd \tau = \frac{B_{j,\sigma}}{k+l}(1-\hC_{j,\sigma}^{k+l}), & &\quad 1 \leq k + l \leq \xi,
\end{alignat}
where
\begin{align}
    C_{i,\tau} = \sum_{j=1}^s \int_0^1 A_{i,\tau,j\sigma}\,\rmd \sigma ,
    \quad
    \hat{C}_{i,\tau} = \sum_{j=1}^s \int_0^1 \hat{A}_{i,\tau,j\sigma}\,\rmd \sigma 
\end{align}
As discussed in~\cite{ch11}, a method satisfying $B(\rho), C(\eta),hC(\eta), D(\xi), \hD(\xi)$ has the order at least $p=\min (\rho,2\eta+2,\xi + \eta + 1)$.

Utilizing simplifying assumptions is effective for deriving energy-preserving CSRK methods with some degrees of freedom for Hamiltonian systems.
However, there is a subtle yet crucial difference between CSRK methods applied to Hamiltonian systems and PCSRK methods applied to Poisson systems.
For Hamiltonian systems, if a CSRK method is energy-preserving, i.e. $M=M^\trans$, then $B(1)$ implies that $B(\rho)$ is satisfied for all $\rho=1,2,\dots$.
In other words, a consistent energy-preserving CSRK method automatically guarantees $B(\rho)$.
This is a very important property, which further simplifies the discussions using simplifying assumptions.
However, this is \emph{not} the case for Poisson systems.
\end{remark}

Now we set $s=3$, and our strategy is as follows.
First, as a sufficient condition for the method to be energy-preserving and of order 4, we require that the method coincides with the 3-degree energy-preserving CSRK methods \eqref{eq:M3alpha}.
We also require the sufficient conditions for being energy-preserving (\cref{thm:pcsrk_energy}) and symmetric (\cref{thm:pcsrk_symmetry}).
The remaining task is to ensure the order conditions for the trees with three vertices and characterize $M_1,M_2,M_3$ and $c_1,c_2,c_3$ satisfying all assumptions and requirements as simply as possible.

As a sufficient condition, we require that 
\begin{equation}
    \label{sufic_cond_M_123}
    M = M_1 + M_2 + M_3 = 
    \begin{bmatrix}
    \tilde{\alpha}+4 & -6\tilde{\alpha} - 6 & 6\tilde{\alpha} \\
    -6\tilde{\alpha} - 6 & 36\tilde{\alpha} + 12 & -36\tilde{\alpha} \\
    6\tilde{\alpha} & -36\tilde{\alpha} & 36\tilde{\alpha}
    \end{bmatrix}.
\end{equation}
In addition to this, we require that $M_1,M_2,M_3$ are symmetric, and the method itself is symmetric.
Under the assumption that $M_1,M_2,M_3$ are symmetric, \cref{thm:pcsrk_symmetry} indicates that the method is symmetric if  
\begin{equation}
    \begin{bmatrix}
    1 & 1 & 1 \\
    0 & -1 & -2 \\
    0 & 0 & 1
    \end{bmatrix}
    M_3
    \begin{bmatrix}
    1 & 0 & 0 \\
    1 & -1 & 0 \\
    1 & -2 & 1
    \end{bmatrix}
    =M_1,
\end{equation}
and $c_2 = 1/2$, $c_1 + c_3 = 1$.
Note that $M_2$ can be arbitrary as long as it is symmetric.

To ensure that the method is of order 4, 
the method needs to satisfy the order conditions for the trees
$\black\tree{\white\vertex1{0.6}\black\vertex2{-0.6}},\black\tree{\white\vertex1{0.6}\white\vertex2{-0.6}},\tree{\vertex10\change\vertex10},
\tree{\change\vertex10\change\vertex10},
\tree{\change\vertex10\vertex10}$.
Note that the order conditions for $\black\tree{\vertex1{0.6}\vertex2{-0.6}}$ and $\tree{\vertex10\vertex10}$ are automatically satisfied because assuming \eqref{sufic_cond_M_123} means that the method is already of order at least 4 when applied to a Hamiltonian system. 
More precisely, the order conditions for the trees of order up to 4 which have only black nodes are automatically satisfied.

The assumptions of the following proposition aid in constructing the intended integrators.

\begin{proposition}
\label{prop:eph}
Assume that a 3-degree PCSRK method satisfies
\begin{align}
& M_i = M_i^\trans, \quad i=1,2,3, \label{method1-1} \\ 
&     {\renewcommand{\arraystretch}{1.2}
\begin{bmatrix}
1 & \tfrac{1}{2} & \tfrac{1}{3} \\
\tfrac12 & \tfrac13 & \tfrac14 \\
\tfrac13 & \tfrac14 & \alpha
\end{bmatrix}}
M = I_3, \quad \alpha \in \bbR, \quad M = M_1+M_2 + M_3,\label{method1-2}\\
&     \begin{bmatrix}
    1 & 1 & 1 \\
    0 & -1 & -2 \\
    0 & 0 & 1
    \end{bmatrix}
    M_3
    \begin{bmatrix}
    1 & 0 & 0 \\
    1 & -1 & 0 \\
    1 & -2 & 1
    \end{bmatrix}
    =M_1,\label{method1-3}\\
& c_1 + c_3 = 1, \quad c_2 = \tfrac12,\label{method1-4}\\
& (c_1 M_1 + c_2 M_2 + c_3 M_3) 
\begin{bmatrix}
1 \\ 1/2 \\ 1/3
\end{bmatrix}
= \begin{bmatrix}
0 \\ 1 \\ 0
\end{bmatrix},\label{method1-5}\\
& 
\begin{bmatrix}
1 & 1/2  & 1/3
\end{bmatrix}
(c_1^2 M_1 + c_2^2 M_2 + c_3^2 M_3) 
\begin{bmatrix}
1 \\ 1/2 \\ 1/3
\end{bmatrix}
= \tfrac{1}{3}.\label{method1-6}
\end{align}
Then, when applied to Poisson systems, the method is energy-preserving, symmetric, and of order at least 4.
\end{proposition}

\begin{proof}
From \cref{thm:pcsrk_energy}, the condition \eqref{method1-1} ensures the energy-preservation.
The condition \eqref{method1-2} guarantees that the method is of order at least 4 when applied to the cases where $S$ is constant, indicating that for Poisson systems, the order conditions for $\black\tree{},
\tree{\vertex10},
\black\tree{\vertex1{0.6}\vertex2{-0.6}},
\tree{\vertex10\vertex10},
\black\tree{\vertex1{1}\vertex2{0}\vertex3{-1}},
\black\tree{\vertex10{\vertex1{0.6}\vertex2{-0.6}}},
\black\tree{\vertex1{-0.6}\vertex2{0.6}\vertex10},
\tree{\vertex10\vertex10,\vertex10}
$
are automatically satisfied.
The conditions \eqref{method1-3} and \eqref{method1-4} relate to the method being symmetric, suggesting that only trees with three vertices must taken into account to guarantee that the method is of order at least 4.

Let $\phi$ be the elementary differentials.
It remains to show that 
$\phi(\black\tree{\white\vertex1{0.6}\black\vertex2{-0.6}})
=
\phi(\black\tree{\white\vertex1{0.6}\white\vertex2{-0.6}}) = 1/3$ and
$
\phi(\tree{\vertex10\change\vertex10}) 
=
\phi(\tree{\change\vertex10\change\vertex10})
=
\phi(\tree{\change\vertex10\vertex10}) 
=1/6
$.
We notice that from \eqref{method1-1} and \eqref{method1-2} that 
$B_{i,\tau} = [1, \ \tau, \ \tau^2] M_i [1, \ 1/2, \ 1/3]^\trans$, $C_{i,\tau} = \tau$.
\begin{itemize}
    \item[\black\tree{\white\vertex1{0.6}\black\vertex2{-0.6}}:] 
    From \eqref{method1-5}, we see that
    \begin{align}
        \phi(\black\tree{\white\vertex1{0.6}\black\vertex2{-0.6}})
        &= \sum_{i=1}^3 \int_0^1 B_{i,\tau} C_{i,\tau} \hC_{i,\tau} \,\rmd \tau \\
        &=
        \int_0^1 \begin{bmatrix}
        \tau & \tau^2 & \tau^3
        \end{bmatrix}
        (c_1M_1 + c_2 M_2 + c_3 M_3) 
        \begin{bmatrix}
        1 \\ 1/2 \\ 1/3
        \end{bmatrix}
        \,\rmd \tau \\
        &= \int_0^1 \begin{bmatrix}
        \tau & \tau^2 & \tau^3
        \end{bmatrix}
        \begin{bmatrix}
        0 \\ 1 \\ 0
        \end{bmatrix}
        \,\rmd \tau
        =\int_0^1 \tau^2 \,\rmd \tau
        =1/3.
    \end{align}
    \item[\black\tree{\white\vertex1{0.6}\white\vertex2{-0.6}}:] From \eqref{method1-6}, we see that
    \begin{align}
        \phi(\black\tree{\white\vertex1{0.6}\white\vertex2{-0.6}})
        &=
        \sum_{i=1}^3 \int_0^1 B_{i,\tau} \hC_{i,\tau}^2\,\rmd \tau \\
        &= 
        \int_0^1 \begin{bmatrix}
        1 & \tau & \tau^2
        \end{bmatrix}
        (c_1^2 M_1 + c_2^2 M_2 + c_3^2 M_3)
        \begin{bmatrix}
        1\\1/2\\ 1/3
        \end{bmatrix}
        \,\rmd \tau \\
        &=
        \begin{bmatrix}
        1 & \tfrac12 & \tfrac13
        \end{bmatrix}
        (c_1^2 M_1 + c_2^2 M_2 + c_3^2 M_3)
        \begin{bmatrix}
        1\\1/2\\ 1/3
        \end{bmatrix}
        =1/3.
    \end{align}
    \item[\tree{\vertex10\change\vertex10}:]
    From \eqref{method1-5} and $\sum_{i=1}^3 B_{i,\tau}=1$, we see that
    \begin{align}
        \phi(\tree{\vertex10\change\vertex10})
        &=
        \sum_{i,j=1}^3 \int_0^1 \int_0^1 B_{i,\tau} A_{i,\tau,j,\sigma} \hC_{j,\sigma}\, \rmd \tau \rmd \sigma \\
        &=
        \int_0^1
        \int_0^1 
        \begin{bmatrix}
        1 & \tfrac{\tau^2}{2} & \tfrac{\tau^3}{3}
        \end{bmatrix}
        (c_1M_1 + c_2 M_2 + c_3 M_3)
        \begin{bmatrix}
        1 \\ 
        \sigma \\ \sigma^2
        \end{bmatrix}
        \,\rmd\tau\rmd \sigma \\
        &=
        \int_0^1 
        \begin{bmatrix}
        1 & \tfrac{\tau^2}{2} & \tfrac{\tau^3}{3}
        \end{bmatrix}
        (c_1M_1 + c_2 M_2 + c_3 M_3)
        \begin{bmatrix}
        1 \\ 
        1/2 \\ 1/3
        \end{bmatrix}
        \,\rmd\tau \\
        &=
        \int_0^1 
        \begin{bmatrix}
        1 & \tfrac{\tau^2}{2} & \tfrac{\tau^3}{3}
        \end{bmatrix}
        \begin{bmatrix}
        0 \\ 
        1 \\ 0
        \end{bmatrix}
        \,\rmd\tau 
        =
        \int_0^1 \tfrac{\tau^2}{2}\,\rmd \tau = \tfrac{1}{6}.
    \end{align} 
    \item[\tree{\change\vertex10\change\vertex10}:]
    From \eqref{method1-2} and \eqref{method1-6}, we see that
    \begin{align}
        \phi(\tree{\change\vertex10\change\vertex10})
        &=
        \sum_{i,j=1}^3
        \int_0^1
        \int_0^1
        B_{i,\tau} \hA_{i,\tau,j,\sigma} C_{j,\sigma}\,\rmd \tau \rmd \sigma \\
        &=
        \sum_{i=1}^3
        \int_0^1
        \int_0^1 
        \begin{bmatrix}
        1 & \tau & \tau^2
        \end{bmatrix}
        M_i
        \begin{bmatrix}
        1 \\ 1/2 \\ 1/3
        \end{bmatrix}
        \begin{bmatrix}
        c_i & \tfrac{c_i^2}{2} & \tfrac{c_i^3}{3}
        \end{bmatrix}
        M
        \begin{bmatrix}
        \sigma \\ \sigma^2 \\ \sigma^3
        \end{bmatrix}
        \,\rmd\tau\rmd\sigma \\
        &=
        \sum_{i=1}^3
        \int_0^1 
        \begin{bmatrix}
        1 & \tau & \tau^2
        \end{bmatrix}
        M_i
        \begin{bmatrix}
        1 \\ 1/2 \\ 1/3
        \end{bmatrix}
        \begin{bmatrix}
        c_i & \tfrac{c_i^2}{2} & \tfrac{c_i^3}{3}
        \end{bmatrix}
        M
        \begin{bmatrix}
        1/2 \\ 1/3 \\ 1/4
        \end{bmatrix}
        \,\rmd\tau \\
        &=
        \sum_{i=1}^3
        \int_0^1 
        \begin{bmatrix}
        1 & \tau & \tau^2
        \end{bmatrix}
        M_i
        \begin{bmatrix}
        1 \\ 1/2 \\ 1/3
        \end{bmatrix}
        \begin{bmatrix}
        c_i & \tfrac{c_i^2}{2} & \tfrac{c_i^3}{3}
        \end{bmatrix}
        \begin{bmatrix}
        0 \\ 1 \\ 0
        \end{bmatrix}
        \,\rmd\tau \\
        &=
        \tfrac{1}{2}
        \int_0^1 
        \begin{bmatrix}
        1 & \tau & \tau^2
        \end{bmatrix}
        (c_1^2 M_1 + c_2^2 M_2 + C_3^2 M_3)
        \begin{bmatrix}
        1 \\ 1/2 \\ 1/3
        \end{bmatrix}
        \,\rmd \tau \\
        &=
        \tfrac{1}{2}
        \int_0^1 
        \begin{bmatrix}
        1 & \tau & \tau^2
        \end{bmatrix}
        \begin{bmatrix}
        0 \\ 0 \\ 1
        \end{bmatrix}
        \,\rmd \tau 
        =
        \tfrac{1}{2}
        \int_0^1 \tau^2 \,\rmd \tau = \tfrac{1}{6}.
    \end{align}
    \item[\tree{\change\vertex10\vertex10}:]
    From \eqref{method1-5} and \eqref{method1-6}, we see that
    \begin{align}
        \phi(\tree{\change\vertex10\vertex10})
        &=
        \sum_{i,j=1}^3
        \int_0^1
        \int_0^1
        B_{i,\tau} \hA_{i,\tau,j,\sigma} \hC_{j,\sigma}\,\rmd \tau \rmd \sigma \\
        &=
        \sum_{i=1}^3
        \int_0^1 
        \int_0^1 
        \begin{bmatrix}
        1 & \tau^2 & \tau^3
        \end{bmatrix}
        M_i
        \begin{bmatrix}
        1 \\ 1/2 \\ 1/3
        \end{bmatrix}
        \begin{bmatrix}
        c_i & \tfrac{c_i^2}{2} & \tfrac{c_i^3}{3}
        \end{bmatrix}
        (c_1M_1 + c_2 M_2 + c_3 M_3)
        \begin{bmatrix}
        1 \\ \sigma \\ \sigma^2
        \end{bmatrix}
        \,\rmd\tau \rmd \sigma\\
        &=
        \sum_{i=1}^3
        \int_0^1 
        \begin{bmatrix}
        1 & \tau^2 & \tau^3
        \end{bmatrix}
        M_i
        \begin{bmatrix}
        1 \\ 1/2 \\ 1/3
        \end{bmatrix}
        \begin{bmatrix}
        c_i & \tfrac{c_i^2}{2} & \tfrac{c_i^3}{3}
        \end{bmatrix}
        (c_1M_1 + c_2 M_2 + c_3 M_3)
        \begin{bmatrix}
        1 \\ 1/2 \\ 1/3
        \end{bmatrix}
        \,\rmd\tau \\
        &=
        \sum_{i=1}^3
        \int_0^1 
        \begin{bmatrix}
        1 & \tau^2 & \tau^3
        \end{bmatrix}
        M_i
        \begin{bmatrix}
        1 \\ 1/2 \\ 1/3
        \end{bmatrix}
        \begin{bmatrix}
        c_i & \tfrac{c_i^2}{2} & \tfrac{c_i^3}{3}
        \end{bmatrix}
        \begin{bmatrix}
        0 \\ 1 \\ 0
        \end{bmatrix}
        \,\rmd\tau = 1/6.
    \end{align} 
\end{itemize}
\end{proof}

The condition \eqref{method1-5} comprises three equations.
However, they are not independent under \eqref{method1-2}, \eqref{method1-3} and \eqref{method1-4}: they actually impose only a single constraint.
We explain this in detail below.

Using $M=M_1+M_2+M_3$ and \eqref{method1-4}, we observe that
\begin{equation}
    c_1M_1 + c_2 M_2 + c_3 M_3 
    = 
    \tfrac{1}{2}M + (c_1 - \tfrac{1}{2})M_1 + (\tfrac{1}{2}-c_1) M_3.
\end{equation}
Applying \eqref{method1-2} and \eqref{method1-3}, which yields
$M[1,\, 1/2,\, 1/3]^\trans = [1,\, 0,\, 0]^\trans$, 
we obtain
\begin{align}
    & (c_1M_1 + c_2 M_2 + c_3 M_3)
    \begin{bmatrix}
    1 \\ 1/2 \\ 1/3
    \end{bmatrix}
    =
    \begin{bmatrix}
    1/2 \\ 0 \\ 0
    \end{bmatrix}
    +
    (c_1 - \tfrac{1}{2})(M_1-M_3)
    \begin{bmatrix}
    1 \\ 1/2 \\ 1/3
    \end{bmatrix} \\
    &\quad =
    \begin{bmatrix}
    1/2 \\ 0 \\ 0
    \end{bmatrix}
    +
    (c_1 - \tfrac{1}{2})
    \begin{bmatrix}
    0 & 1 & 1 \\
    0 & -2 & -2 \\
    0 & 0 & 0
    \end{bmatrix}
    M_3
    \begin{bmatrix}
    1 \\ 1/2 \\ 1/3
    \end{bmatrix}.
\end{align}
Here, we used 
\begin{equation}
    M_1 \begin{bmatrix}
    1 \\ 1/2 \\ 1/3
    \end{bmatrix}
    =
    \begin{bmatrix}
    1 & 1 & 1 \\
    0 & -1 & -2 \\
    0 & 0 & 1
    \end{bmatrix}
    M_3
    \begin{bmatrix}
    1 & 0 & 0 \\
    1 & -1 & 0 \\
    1 & -2 & 1
    \end{bmatrix}
    \begin{bmatrix}
    1 \\ 1/2 \\ 1/3
    \end{bmatrix}
    =
    \begin{bmatrix}
    1 & 1 & 1 \\
    0 & -1 & -2 \\
    0 & 0 & 1
    \end{bmatrix}
    M_3
    \begin{bmatrix}
    1 \\ 1/2 \\ 1/3
    \end{bmatrix}.
\end{equation}
As a result, the condition \eqref{method1-5} can be rewritten as
\begin{align}
    (c_1 - \tfrac{1}{2})
    \begin{bmatrix}
    0 & 1 & 1 \\
    0 & -2 & -2 \\
    0 & 0 & 0
    \end{bmatrix}
    M_3
    \begin{bmatrix}
    1 \\ 1/2 \\ 1/3
    \end{bmatrix}
    =
    \begin{bmatrix}
    -1/2 \\  1 \\ 0
    \end{bmatrix},
\end{align}
and the first and second constraints (rows) are evidently compatible.
Consequently, we only need to consider
\begin{equation}
    \label{method1_M3_cond1}
    (2c_1-1) [0,\, 1 , \, 1] M_3 
    \begin{bmatrix}
    1 \\  1/2 \\ 1/3
    \end{bmatrix} = -1.
\end{equation}

The remaining task is to determine or characterize $M_1, M_2, M_3$ and $c_1,c_2,c_3$ so that they satisfy the conditions \eqref{method1-1}--\eqref{method1-6}. 
The condition \eqref{method1-6} can be rewritten as
\begin{align}
\label{method1_M3_cond2}
(2c_1 -1)^2 \begin{bmatrix}
1 & \tfrac{1}{2} & \tfrac{1}{3}
\end{bmatrix}
M_3
 \begin{bmatrix}
    1 \\  1/2 \\ 1/3
    \end{bmatrix}
    =\tfrac{1}{6}.
\end{align}
It is found that the symmetric matrix $M_3$ satisfying 
\eqref{method1_M3_cond1} and \eqref{method1_M3_cond2}
can be expressed as
\begin{align}
    M_3 &= 
    \begin{bmatrix}
    \tfrac{1}{6(2c_1-1)^2} + \tfrac{1}{2c_1-1} & -\tfrac{1}{2c_1-1} & 0 \\
    -\tfrac{1}{2c_1-1} & 0 & 0\\
    0 & 0 & 0
    \end{bmatrix}
    +
    \gamma_1 
    \begin{bmatrix}
    1 & -3 & 3 \\
    -3 & 0 & 0\\
    3 & 0 & 0
    \end{bmatrix}
    +
    \gamma_2
    \begin{bmatrix}
    1 & -2 & 0 \\
    -2 & 4 & 0\\
    0 & 0 & 0
    \end{bmatrix} \\
    &\phantom{=}\quad
    +
    \gamma_3
    \begin{bmatrix}
    3 & -5 & 0 \\
    -5 & 0 & 6\\
    0 & 6 & 0
    \end{bmatrix}
    +
    \gamma_4
    \begin{bmatrix}
    2 & -3 & 0 \\
    -3 & 0 & 0\\
    0 & 0 & 9
    \end{bmatrix}.
\end{align}
We can set $c_1<1/2$ arbitrary as long as $c_1\neq 0$, and $\gamma_1,\gamma_2,\gamma_3,\gamma_4$ are arbitrary real numbers.
Accordingly, $M_1$ is determined from \eqref{method1-3}.
Thus, we can regard $c_1$, $\gamma_1,\gamma_2,\gamma_3,\gamma_4$ and $\alpha$ in \eqref{method1-2} as free parameters of the fourth-order methods.

\begin{remark}

The exploration of PCSRK methods with a degree greater than three, which are energy-preserving and have an order of at least four, 
is a challenge. Of course,
formulating the conditions corresponding to \eqref{method1-1}--\eqref{method1-6} is relatively straightforward.
However, special attention must be paid to the \eqref{method1-5}-type condition
due to the fact that
\begin{equation}
    \label{method1-5-type}
    (c_1M_1+ \cdots +c_sM_s) \begin{bmatrix}
    1 \\ 1/2 \\1/3 \\ \vdots \\ 1/s
    \end{bmatrix}
    =
    \begin{bmatrix}
    0 \\ 1 \\ 0 \\ \vdots \\ 0
    \end{bmatrix}.
\end{equation}
This equation no longer represents a single constraint.
For example, when $s=4$ or $5$, the condition implies two independent constraints.
This condition corresponds to the simplifying assumption $C(\eta)$ with $k=2$ and $l=1$, which allows for a reduction in the number of order conditions.
Specifically, in the above derivation, $C(\eta)$ with $k=2$ and $l=1$ constitutes a single constraint, which simplifies the calculation.
When $s>3$, it is crucial to carefully consider whether to employ the assumption or treat $\phi(\black\tree{\white\vertex1{0.6}\black\vertex2{-0.6}})$, $\phi(\tree{\vertex10\change\vertex10})$ and $ \phi(\tree{\change\vertex10\vertex10})$ independently.
The advantage of imposing \eqref{method1-5-type}  lies in the fact that the condition for $M_i$ remains linear.

The author believes that a similar approach might help us construct higher-order methods. However, in order to streamline the derivation process, simplifying assumptions or their variants should be incorporated. Some challenges to address include expressing these assumptions in terms of the $M_i$ matrices.

This paper has been focusing on the derivation of 3-degree fourth-order methods due to their apparent practical utility.
\end{remark}

\begin{remark}
The discussion presented above pertains to the case where $C_{i,\tau} = \tau$. However, more general cases can also be examined, as has been discussed for Hamiltonian systems~\cite{mb16}. Further exploration of these cases is not pursued here due to the cumbersome nature of the presentation and the limited practical advantages they appear to offer compared to the methods derived above.
\end{remark}

\begin{remark}
A function $C(y)$ is referred to as a Casimir function if the condition $\nabla C(y)^\trans B(y)=0$ holds for all $y$. In cases where the Casimir takes the quadratic form $C(y) = y^\trans A y$ with a symmetric constant matrix $A$, the $s$-degree $2s$-order PCSRK method~\cite{ch11} exactly inherits the  Casimir. However, the newly introduced fourth-order integrators cannot generally preserve the Casimir. This limitation may be considered a potential drawback of the new family.
\end{remark}

\section{Parameter selection}
\label{sec:4}

This section addresses implementation concerns. We begin by outlining a strategy for implementing PCSRK methods in a general context and then proceed to investigate parallelizable integrators.
As shown below, the parallelizability is characterized in terms of only $\alpha$.
We also investigate the choice of other parameters.

\subsection{Solving the nonlinear equations system}

Let us express $Y_\tau$ as
\begin{equation}
    Y_\tau = y_0 l_0 (\tau) + \sum_{i=1}^s Y_{c_i} l_i (\tau),
\end{equation}
where $l_i (\tau)$ is defined by
\begin{equation}
    l_i (\tau) = \prod_{j=0,\, j \neq i}^s \frac{\tau-c_j}{c_i-c_j}, \quad i = 0,1,\dots,s
\end{equation}
with $c_0 = 0$.
Note that $Y_{c_i} = Z_i$.
We now regard \eqref{pcsrk1} and \eqref{pcsrk2} as a system of nonlinear equations in terms of $Y_{c_1},\dots,Y_{c_s}$:
\begin{equation}
    Y_{c_i} = y_0 + h \sum_{j=1}^s \int_0^1 A_{i,c_i,j,\sigma} S(Y_{c_i}) \nabla H (Y_\sigma)\,\rmd \sigma, \quad i = 1,\dots,s.
\end{equation}
It is not necessary to evaluate \eqref{pcsrk1} at $c_1,\dots,c_s$ as long as it is evaluated at $s$ distinct points; however, $c_1,\dots,c_s$ are considered to avoid cumbersome presentation.
Let
\begin{equation}
    Y = \begin{bmatrix}
    Y_{c_1} \\
    \vdots \\
    Y_{c_s}
    \end{bmatrix}
    \in \bbR^{sN}.
\end{equation}
Then, we are concerned with solving
\begin{equation}
    \Phi(Y) = 
    Y - e_s \otimes y_0 -
    h
    \begin{bmatrix}\displaystyle
    \sum_{j=1}^s \int_0^1 A_{1,c_1,j,\sigma} S(Y_{c_i}) \nabla H (Y_\sigma) \,\rmd \sigma \\
    \vdots\\
    \displaystyle
    \sum_{j=1}^s \int_0^1 A_{s,c_s,j,\sigma} S(Y_{c_i}) \nabla H (Y_\sigma) \,\rmd \sigma
    \end{bmatrix}
    = 0.
\end{equation}
Note here that
\begin{equation}
    \nabla_{Y_{c_j}} \big( S(Y_{c_i}) \nabla H(Y_\sigma)\big)
    \approx
    S(y_0) \nabla^2 H (Y_\sigma) l_j (\sigma)
    \approx
    S(y_0) \nabla^2 H (y_0) l_j (\sigma).
\end{equation}
In the first `$\approx$', the term involving $\nabla _{Y_{c_j}} S(Y_{c_i})$ is omitted to avoid that tensors appear,
by taking in mind that even if we care about such a term, the following discussion remains true.
Thus, we have
\begin{align}
    \Phi ' (Y) &\approx  I_{sN} - h
    \begin{bmatrix}
    \sum_{j=1}^s\int_0^1 A_{1,c_1,j,\sigma} S(y_0) \nabla^2 H(y_0) l_1(\sigma) 
    \,\rmd \sigma& \cdots & 
    \sum_{j=1}^s\int_0^1 A_{1,c_1,j,\sigma} S(y_0) \nabla^2 H(y_0) l_s(\sigma)\,\rmd \sigma \\
    \vdots & \ddots & \vdots \\
    \sum_{j=1}^s\int_0^1 A_{s,c_s,j,\sigma} S(y_0) \nabla^2 H(y_0) l_1(\sigma)\,\rmd \sigma & \cdots & 
    \sum_{j=1}^s\int_0^1 A_{s,c_s,j,\sigma} S(y_0) \nabla^2 H(y_0) l_s(\sigma)\,\rmd \sigma
    \end{bmatrix} \\
    &=
    I_{sN} - h E \otimes J_0,
\end{align}
where
\begin{equation}
    \label{eq:defE}
    E_{ij} = \sum_{k=1}^s \int_0^1 A_{1,c_i,k,\sigma} l_j (\sigma) \,\rmd \sigma
    = \int_0^1 \big( \sum_{k=1}^s A_{1,c_i,k,\sigma}\big) l_j (\sigma)\,\rmd \sigma
\end{equation}
and 
\begin{equation}
    J_0 = S(y_0) \nabla^2 H(y_0) 
\end{equation}
denotes an approximate Jacobian matrix.
Therefore, a simplified Newton-like method gives the iteration formula
\begin{equation}
    \label{eq:sNewton}
    (I_{sN} - h E\otimes J_0) \rho^l = - \Phi (Y^l), 
    \quad Y^{l+1} = Y^l + \rho^l, \quad l = 0,1,2,\dots.
\end{equation}

\subsection{Parallelizable integrators}

As is the case with implicit Runge--Kutta methods and continuous stage Runge--Kutta methods,
if the matrix $E$ has only real, distinct eigenvalues,
the linear system \eqref{eq:sNewton} of size $sN$ can be computed efficiently by using parallel architectures.
If all eigenvalues of $E$ are real and distinct, there exists a matrix $T$ such that 
\begin{equation}
    T^{-1} E T = \diag (\lambda_1,\dots,\lambda_s), \quad 
    \lambda_1,\dots,\lambda_s\in \bbR.
\end{equation}
Let $Q = (I_{sN} - h E \otimes J_0)$ and define
$\overline{Q} := (T^{-1}\otimes I_N) Q (T\otimes I_N)$.
Then, it follows that 
\begin{equation}
    \overline{Q} = \diag(I_N - h\lambda_1 J_0, \dots, I_N - h \lambda_s J_0),
\end{equation}
which is block diagonal.
Therefore, $\rho^l$ can be calculated based on
\begin{align}
    \overline{\Phi} (Y_l) &= (T^{-1} \otimes I_N) \Phi (Y^l), \\
    \overline{Q} \overline{\rho}^l &= - \overline{\Phi} (Y^l),\\
    \overline{\rho}^l &= (T^{-1} \otimes I_N) \rho ^l.
\end{align}
The key point is that the linear system for $\overline{\rho}^l$ of size $sN$ consists of $s$ linear systems of size $N$.
Thus, the most computationally heavy part in updating $Y_l$ can be computed in parallel.
We note that $T$ is a $3$-by-$3$ matrix and thus the explicit computation of $T^{-1}$ and the multiplication of $T^{-1}\otimes I_N$ with a vector are not particularly challenging.

We note that $E$ is identical to that appearing in the study of energy-preserving methods for  Hamiltonian system~\cite{mb16}.

\begin{theorem}[\cite{mb16}]
The eigenvalues of the matrix $E$ defined in \eqref{eq:defE} are independent of the $c_i$ values.
\end{theorem}

This theorem indicates that the parallelizability does not depend on the $c_i$ values.
If the method is parallelizable when applied to Hamiltonian systems, it is also parallelizable to Poisson systems.

For the method satisfying the assumptions in \cref{prop:eph}, according to \cite[Section~5.3.1]{mb16},
the corresponding $E$ has real and distinct eigenvalues if 
\begin{equation}
    \label{eq:cond_real}
    -\frac{\tilde{\alpha}}{300}
    >
    \frac{1}{6} 2^{2/3}
    + \frac{5}{24} 2^{1/3} + \frac{1}{4}
    \approx
    0.7770503941
    \quad \text{with }
    \tilde{\alpha} = \frac{1}{36\alpha-7}.
\end{equation}

\subsection{Discussion for the accuracy}
\label{subsec:accuracy}

For parallelizable fourth-order integrators, the most computationally demanding aspect is solving linear systems of size $d$. In contrast,  for the energy-preserving 2-degree PCSRK method, one needs to solve a linear system of size $2d$.
When dealing with dense coefficient matrices,
direct solvers are often preferred for solving such linear systems.
The computational complexity of direct solvers is proportional to the cube of the size of the linear systems.

Considering this,
we assume that our parallelizable integrators are eight times faster than the 2-degree PCSRK method. However, we note that this is not always the case, as actual computational costs can vary significantly based on the specific problem and the choice of linear solver.

We compare the coefficients of elementary differentials with trees of order 5 for the P-series expansions of the exact solution, 2-degree PCSRK method (AVF(4)), and the proposed method.
Table~\ref{table:black} compares the coefficients for the trees for which every vertex is black.
Clearly, for both AVF(4) and the proposed methods, the coefficients for the trees \black\tree{\vertex1{-0.6}\vertex2{0.6}{\vertex1{0.6}\vertex2{-0.6}}},
\black\tree{\vertex1{-0.6}\vertex2{0.6}\vertex10\vertex10},
\black\tree{\vertex10\vertex10\vertex1{-0.6}\vertex2{0.6}},
\black\tree{\vertex10\vertex10\vertex10\vertex10}
differ from those for the exact solution.

However, the proposed method exhibits an interesting property for the coefficients for the bi-coloured tree.
Tables~\ref{table:9-1}--\ref{table:9-9} show the coefficients for bi-coloured trees with black roots with the following quantities:
\begin{align}
    &\frac{(2c_1-1)^2 (\gamma_1+\gamma_3+\gamma_4)}{120},
    \label{A}\tag{A} \\
    &\frac{c_1^2}{12} - \frac{c_1}{12} + \frac{5}{24},
    \label{B}\tag{B} \\
    &\frac{5}{72} - \frac{\tilde{\alpha}}{1800},
    \label{C}\tag{C}\\
    &\frac{c_1(c_1-1)\tilde{\alpha}}{180},
    \label{D}\tag{D}
    \\
    &\frac{(2c_1-1)(2\gamma_3+3\gamma_4)}{1440},
    \label{E}\tag{E}
    \\
    &\frac{c_1(c_1-1)(2c_1-1)(2\gamma_3+3\gamma_4)}{144},
    \label{F}\tag{F}
    \\
    &\frac{c_1(c_1-1)(2c_1-1)^2(\gamma_1+\gamma_3 + \gamma_4)}{24},
    \label{G}\tag{G}\\
    &\frac{(2c_1-1)^2(4\gamma_2+12\gamma_3 + 9\gamma_4)}{288},
    \label{H}\tag{H}
\end{align}
We note that for AVF(4), even for the trees of the form \black\tree{\vertex1{1.2}\vertex2{0.4}\vertex3{-0.4}\vertex4{-1.2}},
\black\tree{\vertex1{-1}\vertex2{0}\vertex3{1}\vertex10},
\tree{\vertex1{0.6}\vertex10\vertex3{-0.6}\vertex10},
\black\tree{\vertex10\vertex1{-1}\vertex2{0}\vertex3{1}},
\black\tree{\vertex10\vertex1{-0.6}\vertex2{0.6}\vertex10}, 
some coefficients for the trees with white nodes differ from the exact ones.
In contrast, for the proposed method, if the parameters satisfy
\begin{equation}
    \label{eq:ABEFGH}
    \eqref{A} = \frac{1}{180},
    \quad 
    \eqref{B} = \frac{1}{5},
    \quad
    \eqref{E} = -\eqref{F} = \eqref{G}  = -\frac{1}{360},
    \quad 
    \eqref{H} = \frac{1}{80},
\end{equation}
all coefficients are exact.
These conditions are compatible by choosing the parameters satisfying
\begin{align}
    & c_1 = \frac{1}{2} - \frac{\sqrt{15}}{10}, \\
    & \gamma_1 + \gamma_3 + \gamma_4 = \frac{10}{9}, \\
    & 2\gamma_3 + 3 \gamma_4 = \frac{20}{\sqrt{15}}, \\
    & 4\gamma_2 + 12\gamma_3 + 9 \gamma_4 = 6.
\end{align}
For the proposed method, only the coefficients for the bi-coloured trees of the shape \black\tree{\vertex1{-0.6}\vertex2{0.6}{\vertex1{0.6}\vertex2{-0.6}}},
\black\tree{\vertex1{-0.6}\vertex2{0.6}\vertex10\vertex10},
\black\tree{\vertex10\vertex10\vertex1{-0.6}\vertex2{0.6}},
\black\tree{\vertex10\vertex10\vertex10\vertex10} that are expressed in terms of \eqref{C} or \eqref{D} cannot coincide with the exact ones if $\tilde{\alpha}$ is in the range of \eqref{eq:cond_real}.

We also note that, apart from constructing efficient fourth-order integrators, if we set $\tilde{\alpha}=5$ all the coefficients for the trees with 5 vertices can be exact.
However, the condition \eqref{eq:ABEFGH} indicates that there remains one degree of freedom.
Therefore, 3-degree sixth-order integrators are not unique, which suggests that there exists a family of $s$-degree $2s$-order energy-preserving integrators (for $s=1$ and $2$, such integrators are unique).

\subsection{Numerical verification}

As a simple test confirming the order of accuracy and the effect of the choice of parameters for the proposed method,
we employ the Lotka--Volterra system
\begin{align}
    & S(y) = 
    \begin{bmatrix}
    0 & c y_1 y_2 & b c y_1 y_3 \\
    -c y_1 y_2 & 0 & - y_2 y_3 \\
    -bcy_1 y_3 & y_2 y_3 & 0
    \end{bmatrix}, \\
    & H(y) = aby_1 + y_2 - a y_3 + \nu \log y_2 - \mu \log y_3.
\end{align}
For the numerical experiment, the parameters were set to $a=-2$, $b=-1$, $c=-0.5$, $\nu=1$, $\mu = 2$, and the initial values $y_0 = (1.0,1.9,0.5)$.
In the numerical experiment, a high-order quadrature with the tolerance $10^{-12}$ is used to approximate the integrals appearing in the integrator.

The integrator actually preserves the Hamiltonian $H(y)$.
For instance, when the step size is set to $h=0.05$, the error between true and numerical Hamiltonians grows as the time integration proceeds due to the rounding errors and the use of quadrature; however, at $t=10$, this error remains below $10^{-12}$.
This level of accuracy is noteworthy, actually highlighting the energy-preservation, especially when compared to the error in the Casimir invariant at the same point in time, which exceeds $0.01$.

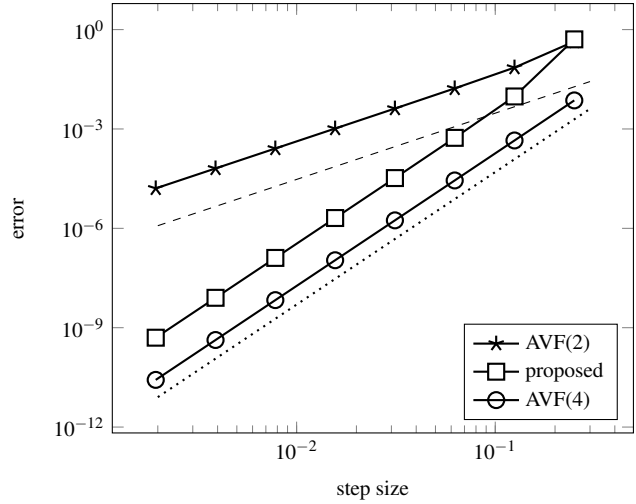
\begin{figure}[t]
    \centering
    \begin{tikzpicture}
    \begin{loglogaxis}[xlabel={step size},ylabel={error},
    legend entries={AVF(2),proposed,AVF(4)},
legend style={legend pos=south east,legend cell align=left},]
    \addplot[thick,mark=star, mark options={fill=white},mark size=3pt,
    ] coordinates {
(0.25, 0.4387107766400442)
(0.125, 0.07138686990743914)
(0.0625, 0.016783848157749297)
(0.03125, 0.004137716809511271)
(0.015625, 0.0010308991239590336)
(0.0078125, 0.0002575058093531552)
(0.00390625, 6.43627941073113e-5)
(0.001953125, 1.6089845222822736e-5)
};
    \addplot[thick,mark=square*, mark options={fill=white},mark size=3pt] coordinates {
(0.25, 0.5077570786574491)
(0.125, 0.009592748765109195)
(0.0625, 0.0005391735621606368)
(0.03125, 3.294929505801119e-5)
(0.015625, 2.048233729961848e-6)
(0.0078125, 1.2785679620238103e-7)
(0.00390625, 7.986870770348962e-9)
(0.001953125, 4.989964260661792e-10)
};
    \addplot[thick,mark=o, mark options={fill=white},mark size=3pt] coordinates {
(0.25, 0.007284847717734343)
(0.125, 0.00044945749653701604)
(0.0625, 2.7908705190459515e-5)
(0.03125, 1.7410912957927447e-6)
(0.015625, 1.0880277696856638e-7)
(0.0078125, 6.8009200251500015e-9)
(0.00390625, 4.248765740052504e-10)
(0.001953125, 2.656180909990692e-11)
};
    \addplot[black,dashed,domain=2e-3:0.3] {0.3*x^2};
    \addplot[black,dotted,thick,domain=2e-3:0.3] {0.5*x^4};
    \end{loglogaxis}
    \end{tikzpicture}
    \caption{Error at $t=1$ of numerical solutions for the Lotka--Volterra system.
    The proposed method with the parameters \eqref{eq:opt_param} and $\tilde{\alpha} = -234$ is compared with the 2nd and 4th order methods proposed in~\cite{ch11}.
    Dashed and dotted lines show the slope for the 2nd and 4th-order convergence, respectively.}
    \label{fig1}
\end{figure}

We now check that the proposed integrator is actually of order 4.
We set $\gamma_1,\dots,\gamma_4$ to
\begin{align}
    (\gamma_1,\gamma_2,\gamma_3,\gamma_4)
    =
    \Big( \frac{10}{3} - \frac{2\sqrt{15}}{3}, \frac{23}{2} - 2\sqrt{15}, -\frac{20}{3} + \frac{2\sqrt{15}}{3}, \frac{40}{9} \Big).
    \label{eq:opt_param}
\end{align}
This choice corresponds to the 6th order integrator presented in~\cite{ch11} when $\tilde{\alpha}=5$.
In Fig.~\ref{fig1}, the error behaviour of the proposed method with this choice of parameters and $\tilde{\alpha} = -234$ are compared with the 2nd and 4th order integrators proposed in~\cite{ch11}.
This figure supports that the proposed integrator actually attains the 4th-order accuracy.
While the error of the proposed method is substantially bigger than that of AVF(4), efficiency could still favour the proposed method if its computation is more than twice as fast as AVF(4).

\begin{remark}
    We discuss the efficiency of our proposed method in more detail. 
    It has been observed that to achieve a comparable level of accuracy, the proposed method necessitates nearly halving the step size. 
    This suggests that if the computation time per time step of the proposed method is at least twice as fast as that of AVF(4), then the proposed method is considered to be more efficient.

    A comparative analysis of the two methods is provided. 
    We examine a scenario where the linear system of size $2d$, required for solving AVF(4), is solved using a direct method on a single core. 
    Conversely, three linear systems of size $d$ for the proposed method are solved using the same direct method but potentially by three cores in parallel. 
    Assuming that the simplified Newton iterations for both methods demand a nearly identical number of iterations to the convergence, the proposed method is approximately eight times faster because the computational cost of a direct solver typically scales with the cube of the system size.

    However, it is important to note that such estimations are contingent upon the specifics of the problem and the selection of linear solvers. 
    In the case of AVF(4), further parallelization might be achievable. Nonetheless, most techniques for parallelization can also be applied to the proposed method. Therefore, expecting the proposed method to compute a single step at least twice as fast as AVF(4) appears to be a reasonable assumption in most instances, particularly for large-scale problems.
\end{remark}

\section{Concluding remarks}
\label{sec:5}

In this paper, we have shown a family of fourth-order energy-preserving integrators based on the partitioned CSRK methods.
The integrators can be implemented in a parallel architecture if one of the parameters is chosen such that it satisfies a certain inequality.
Actual implementation must tailored to the problem at hand and it is anticipated that these integrators, with their parallelizability, will demonstrate substantial efficiency when applied to large-scale problems. 
A more detailed evaluation of their performance in large-scale contexts will be the subject of forthcoming research,
as similar evaluation was studied for the CSRK methods~\cite{sk21}.


\bibliographystyle{spmpsci} 
\bibliography{references}

\appendix
\section{Coefficients of the elementary differentials}

We present the coefficients of the elementary differentials, which were used in \cref{subsec:accuracy}.
The Julia code employed for verifying these coefficients is available on GitHub. 
This code can be accessed at the following repository: \url{https://github.com/yutomiyatake/EP_fourth_Poisson}.

\begin{table}
\caption{Coefficients of the elementary differentials with trees of order 5 for the P-series expansions of the exact solution, the fourth order AVF collocation method, and the proposed method.}
\label{table:black}
{\renewcommand\arraystretch{1.5}
\begin{tabular}{cccccccccc}
     $t$  & \black\tree{\vertex1{1.2}\vertex2{0.4}\vertex3{-0.4}\vertex4{-1.2}}
     & \black\tree{\vertex1{-1}\vertex2{0}\vertex3{1}\vertex10}
     & \black\tree{\vertex1{-0.6}\vertex2{0.6}{\vertex1{0.6}\vertex2{-0.6}}}
     & \black\tree{\vertex1{-0.6}\vertex2{0.6}\vertex10\vertex10}
     & \tree{\vertex1{0.6}\vertex10\vertex3{-0.6}\vertex10}
     & \black\tree{\vertex10\vertex1{-1}\vertex2{0}\vertex3{1}}
     & \black\tree{\vertex10\vertex1{-0.6}\vertex2{0.6}\vertex10}
     & \black\tree{\vertex10\vertex10\vertex1{-0.6}\vertex2{0.6}}
     & \black\tree{\vertex10\vertex10\vertex10\vertex10} \\
     exact solution & $\tfrac15$ & $\tfrac{1}{10}$ & $\tfrac{1}{15}$ & $\tfrac{1}{30}$ & $\tfrac{1}{20}$ & $\tfrac{1}{20}$ & $\tfrac{1}{40}$ & $\tfrac{1}{60}$ & $\tfrac{1}{120}$\\
     \hline
     AVF(4) & $\tfrac15$ & $\tfrac{1}{10}$ & $\tfrac{5}{72}$ & $\tfrac{5}{144}$ & $\tfrac{1}{20}$ & $\tfrac{1}{20}$ & $\tfrac{1}{40}$ & $\tfrac{1}{72}$ & $\tfrac{1}{144}$\\
     difference & $0$ & $0$ & $\tfrac{1}{360}$ & $\tfrac{1}{720}$ & $0$ & $0$ & $0$ & $-\tfrac{1}{360}$ & $-\tfrac{1}{720}$\\
     \hline
     proposed & $\tfrac15$ & $\tfrac{1}{10}$ & $\tfrac{120\theta+5}{72}$ & $\tfrac{12\theta +5}{144}$ & $\tfrac{1}{20}$ & $\tfrac{1}{20}$ & $\tfrac{1}{40}$ & $\tfrac{1-12\theta}{72}$ & $\tfrac{1-12\theta}{144}$\\
     difference & $0$ & $0$ & $\tfrac{60\theta+1}{360}$ & $\tfrac{60\theta+1}{720}$ & $0$ & $0$ & $0$ & $-\tfrac{60\theta+1}{360}$ & $-\tfrac{60\theta+1}{720}$
\end{tabular}
}
\end{table}

\begin{table}
\newsavebox{\tone}
\sbox{\tone}{
    \black\tree{\vertex1{1.2}\vertex2{0.4}\vertex3{-0.4}\vertex4{-1.2}}
}
\caption{Coefficients of the elementary differentials with trees of the shape \usebox{\tone} with the black root for the P-series expansions of the exact solution, the fourth order AVF collocation method, and the proposed method.}
\label{table:9-1}
{\renewcommand\arraystretch{1.5}
\begin{tabular}{cccccc}
$t$
&
\black\tree{\vertex1{1.2}\vertex2{0.4}\vertex3{-0.4}\vertex4{-1.2}}
&
\black\tree{\white\vertex1{1.2}\black\vertex2{0.4}\vertex3{-0.4}\vertex4{-1.2}}
&
\black\tree{\white\vertex1{1.2}\vertex2{0.4}\black\vertex3{-0.4}\vertex4{-1.2}}
&
\black\tree{\white\vertex1{1.2}\vertex2{0.4}\vertex3{-0.4}\black\vertex4{-1.2}}
&
\black\tree{\white\vertex1{1.2}\vertex2{0.4}\vertex3{-0.4}\vertex4{-1.2}} \\
exact solution
&
$\tfrac{1}{5}$
&
$\tfrac{1}{5}$
&
$\tfrac{1}{5}$
&
$\tfrac{1}{5}$
&
$\tfrac{1}{5}$
\\
\hline
AVF(4)
&
$\tfrac{1}{5}$
&
$\tfrac{7}{36}$
&
$\tfrac{7}{36}$
&
$\tfrac{7}{36}$
&
$\tfrac{7}{36}$
\\
\hline
proposed
&
$\tfrac{1}{5}$
&
$\tfrac{1}{5}$
&
\eqref{A} + $\tfrac{7}{36}$
&
\eqref{B}
&
\eqref{B}
\end{tabular}
}
\end{table}

\begin{table}
\newsavebox{\ttwo}
\sbox{\ttwo}{
    \black\tree{\vertex1{-1}\vertex2{0}\vertex3{1}\vertex10}
}
\caption{Coefficients of the elementary differentials with trees of the shape \usebox{\ttwo} with the black root for the P-series expansions of the exact solution, the fourth order AVF collocation method, and the proposed method.}
\label{table:9-2}
{\renewcommand\arraystretch{1.5}
\begin{tabular}{ccccccccccccc}
$t$
&
\black\tree{\vertex1{-1}\vertex2{0}\vertex3{1}\vertex10}
&
\black\tree{\vertex1{-1}\white\vertex2{0}\black\vertex3{1}\vertex10}
&
\black\tree{\vertex1{-1}\vertex2{0}\white\vertex3{1}\black\vertex10}
&
\black\tree{\vertex1{-1}\vertex2{0}\vertex3{1}\white\vertex10}
&
\black\tree{\white\vertex1{-1}\vertex2{0}\black\vertex3{1}\vertex10}
&
\black\tree{\vertex1{-1}\white\vertex2{0}\vertex3{1}\black\vertex10}
&
\black\tree{\vertex1{-1}\white\vertex2{0}\black\vertex3{1}\white\vertex10}
&
\black\tree{\vertex1{-1}\vertex2{0}\white\vertex3{1}\vertex10}
&
\black\tree{\white\vertex1{-1}\black\vertex2{0}\white\vertex3{1}\vertex10}
&
\black\tree{\white\vertex1{-1}\vertex2{0}\black\vertex3{1}\white\vertex10} 
&
\black\tree{\white\vertex1{-1}\vertex2{0}\vertex3{1}\black\vertex10}
&
\black\tree{\white\vertex1{-1}\vertex2{0}\vertex3{1}\vertex10}\\
exact solution
&$\tfrac{1}{10}$&$\tfrac{1}{10}$&$\tfrac{1}{10}$&$\tfrac{1}{10}$&$\tfrac{1}{10}$&$\tfrac{1}{10}$&$\tfrac{1}{10}$&$\tfrac{1}{10}$&$\tfrac{1}{10}$&$\tfrac{1}{10}$&$\tfrac{1}{10}$&$\tfrac{1}{10}$ \\
\hline
AVF(4)
& 
$\tfrac{1}{10}$
& 
$\tfrac{1}{10}$
& 
$\tfrac{7}{72}$
& 
$\tfrac{1}{10}$
& 
$\tfrac{7}{72}$
& 
$\tfrac{7}{72}$
& 
$\tfrac{1}{10}$
& 
$\tfrac{7}{72}$
& 
$\tfrac{7}{72}$
& 
$\tfrac{7}{72}$
& 
$\tfrac{7}{72}$
& 
$\tfrac{7}{72}$\\
\hline
proposed 
&
$\tfrac{1}{10}$
&
$\tfrac{1}{10}$
&
$\frac{\eqref{A}}{2} + \tfrac{7}{72}$
&
$\tfrac{1}{10}$
&
$\frac{\eqref{A}}{2} + \tfrac{7}{72}$
&
$\frac{\eqref{B}}{2}$
&
$\tfrac{1}{10}$
&
$\frac{\eqref{A}}{2} + \tfrac{7}{72}$
&
$\frac{\eqref{B}}{2}$
&
$\frac{\eqref{A}}{2} + \tfrac{7}{72}$
&
$\frac{\eqref{B}}{2}$
&
$\frac{\eqref{B}}{2}$

\end{tabular}
}
\end{table}

\begin{table}
\newsavebox{\tthree}
\sbox{\tthree}{
    \black\tree{\vertex1{-0.6}\vertex2{0.6}{\vertex1{0.6}\vertex2{-0.6}}}
}
\caption{Coefficients of the elementary differentials with trees of the shape \usebox{\tthree} with the black root for the P-series expansions of the exact solution, the fourth order AVF collocation method, and the proposed method.}
\label{table:9-3}

{\renewcommand\arraystretch{1.5}
\begin{tabular}{cccccccccccc}
$t$
&
\black\tree{\vertex1{-0.6}\vertex2{0.6}{\vertex1{0.6}\vertex2{-0.6}}}
&
\black\tree{\white\vertex1{-0.6}\black\vertex2{0.6}{\vertex1{0.6}\vertex2{-0.6}}}
&
\black\tree{\vertex1{-0.6}\white\vertex2{0.6}{\black\vertex1{0.6}\vertex2{-0.6}}}
&
\black\tree{\vertex1{-0.6}\vertex2{0.6}{\white\vertex1{0.6}\black\vertex2{-0.6}}}
&
\black\tree{\white\vertex1{-0.6}\vertex2{0.6}{\black\vertex1{0.6}\vertex2{-0.6}}}
&
\black\tree{\white\vertex1{-0.6}\black\vertex2{0.6}{\white\vertex1{0.6}\black\vertex2{-0.6}}}
&
\black\tree{\vertex1{-0.6}\white\vertex2{0.6}{\vertex1{0.6}\black\vertex2{-0.6}}}
&
\black\tree{\vertex1{-0.6}\vertex2{0.6}{\white\vertex1{0.6}\vertex2{-0.6}}}
&
\black\tree{\vertex1{-0.6}\white\vertex2{0.6}{\vertex1{0.6}\vertex2{-0.6}}}
& \\
exact solution
&$\tfrac{1}{15}$&$\tfrac{1}{15}$&$\tfrac{1}{15}$&$\tfrac{1}{15}$&$\tfrac{1}{15}$&$\tfrac{1}{15}$&$\tfrac{1}{15}$&$\tfrac{1}{15}$&$\tfrac{1}{15}$\\
\hline
AVF(4)
& 
$\tfrac{5}{72}$
& 
$\tfrac{5}{72}$
& 
$\tfrac{5}{72}$
& 
$\tfrac{5}{72}$
& 
$\tfrac{5}{72}$
& 
$\tfrac{5}{72}$
& 
$\tfrac{5}{72}$
& 
$\tfrac{5}{72}$
& 
$\tfrac{5}{72}$
\\
\hline
proposed 
&
\eqref{C}
&
\eqref{C}
&
$\eqref{D}+\tfrac{5}{72}$
&
$\eqref{E}+\tfrac{5}{72}$
&
$\eqref{D}+\tfrac{5}{72}$
&
$\eqref{E}+\tfrac{5}{72}$
&
$-\eqref{F}+\tfrac{5}{72}$
&
$-\frac{\eqref{A}}{2}+\tfrac{5}{72}$
&
$\eqref{G}+\tfrac{5}{72}$
\end{tabular}
}

{\renewcommand\arraystretch{1.5}
\begin{tabular}{ccccccccccccc}
$t$
&
\black\tree{\white\vertex1{-0.6}\black\vertex2{0.6}{\white\vertex1{0.6}\vertex2{-0.6}}}
&
\black\tree{\white\vertex1{-0.6}\vertex2{0.6}{\black\vertex1{0.6}\white\vertex2{-0.6}}}
&
\black\tree{\white\vertex1{-0.6}\vertex2{0.6}{\vertex1{0.6}\vertex2{-0.6}}}\\
exact solution
&$\tfrac{1}{15}$&$\tfrac{1}{15}$&$\tfrac{1}{15}$ \\
\hline
AVF(4)
& 
$\tfrac{5}{72}$
& 
$\tfrac{5}{72}$
& 
$\tfrac{5}{72}$\\
\hline
proposed 
&
$-\frac{\eqref{A}}{2}+\tfrac{5}{72}$
&
$-\eqref{F}+\tfrac{5}{72}$
&
$\eqref{G}+\tfrac{5}{72}$

\end{tabular}
}
\end{table}

\begin{table}
\newsavebox{\tfour}
\sbox{\tfour}{
    \black\tree{\vertex1{-0.6}\vertex2{0.6}\vertex10\vertex10}
}
\caption{Coefficients of the elementary differentials with trees of the shape \usebox{\tfour} with the black root for the P-series expansions of the exact solution, the fourth order AVF collocation method, and the proposed method.}
\label{table:9-4}
{\renewcommand\arraystretch{1.5}
\begin{tabular}{ccccccccccccccccc}
$t$
& \black\tree{\vertex1{-0.6}\vertex2{0.6}\vertex10\vertex10}
& \black\tree{\white\vertex1{-0.6}\black\vertex2{0.6}\vertex10\vertex10}
& \black\tree{\vertex1{-0.6}\white\vertex2{0.6}\black\vertex10\vertex10}
& \black\tree{\vertex1{-0.6}\vertex2{0.6}\white\vertex10\black\vertex10}
& \black\tree{\vertex1{-0.6}\vertex2{0.6}\vertex10\white\vertex10}
& \black\tree{\white\vertex1{-0.6}\vertex2{0.6}\black\vertex10\vertex10}
& \black\tree{\white\vertex1{-0.6}\black\vertex2{0.6}\white\vertex10\black\vertex10}
& \black\tree{\white\vertex1{-0.6}\black\vertex2{0.6}\vertex10\white\vertex10}
& \black\tree{\vertex1{-0.6}\white\vertex2{0.6}\vertex10\black\vertex10}
& \black\tree{\vertex1{-0.6}\white\vertex2{0.6}\black\vertex10\white\vertex10}
 \\
exact solution &
$\tfrac{1}{30}$ &
$\tfrac{1}{30}$ &
$\tfrac{1}{30}$ &
$\tfrac{1}{30}$ &
$\tfrac{1}{30}$ &
$\tfrac{1}{30}$ &
$\tfrac{1}{30}$ &
$\tfrac{1}{30}$ &
$\tfrac{1}{30}$ &
$\tfrac{1}{30}$ \\
\hline
AVF(4) &
$\tfrac{5}{144}$ &
$\tfrac{5}{144}$ &
$\tfrac{5}{144}$ &
$\tfrac{5}{144}$ &
$\tfrac{5}{144}$ &
$\tfrac{5}{144}$ &
$\tfrac{5}{144}$ &
$\tfrac{5}{144}$ &
$\tfrac{5}{144}$ &
$\tfrac{5}{144}$ \\
\hline
proposed &
$\tfrac{\eqref{C}}{2}$ &
$\tfrac{\eqref{C}}{2}$ &
$\tfrac{\eqref{D}}{2} + \tfrac{5}{144}$ &
$\tfrac{-\eqref{A}}{4} + \tfrac{5}{144}$ &
$\tfrac{\eqref{C}}{2}$ &
$\tfrac{\eqref{D}}{2} + \tfrac{5}{144}$ &
$\tfrac{-\eqref{A}}{4} + \tfrac{5}{144}$ &
$\tfrac{\eqref{C}}{2}$ &
$\tfrac{\eqref{G}}{2}+\tfrac{5}{144}$ &
$\tfrac{\eqref{D}}{2} + \tfrac{5}{144}$ 

\end{tabular}
}

{\renewcommand\arraystretch{1.5}
\begin{tabular}{ccccccccccccccccc}
$t$
& \black\tree{\vertex1{-0.6}\vertex2{0.6}\white\vertex10\vertex10}
& \black\tree{\vertex1{-0.6}\white\vertex2{0.6}\vertex10\vertex10}
& \black\tree{\white\vertex1{-0.6}\black\vertex2{0.6}\white\vertex10\vertex10}
& \black\tree{\white\vertex1{-0.6}\vertex2{0.6}\black\vertex10\white\vertex10}
& \black\tree{\white\vertex1{-0.6}\vertex2{0.6}\vertex10\black\vertex10}
& \black\tree{\white\vertex1{-0.6}\vertex2{0.6}\vertex10\vertex10} \\
exact solution &
$\tfrac{1}{30}$ &
$\tfrac{1}{30}$ &
$\tfrac{1}{30}$ &
$\tfrac{1}{30}$ &
$\tfrac{1}{30}$ &
$\tfrac{1}{30}$ \\
\hline
AVF(4) &
$\tfrac{5}{144}$ &
$\tfrac{5}{144}$ &
$\tfrac{5}{144}$ &
$\tfrac{5}{144}$ &
$\tfrac{5}{144}$ &
$\tfrac{5}{144}$ \\
\hline
proposed &
$\tfrac{-\eqref{A}}{4} + \tfrac{5}{144}$ &
$\tfrac{\eqref{G}}{2}+\tfrac{5}{144}$ &
$\tfrac{-\eqref{A}}{4} + \tfrac{5}{144}$ &
$\tfrac{\eqref{D}}{2} + \tfrac{5}{144}$ &
$\tfrac{\eqref{G}}{2}+\tfrac{5}{144}$ &
$\tfrac{\eqref{G}}{2}+\tfrac{5}{144}$ 

\end{tabular}
}
\end{table}

\begin{table}
\newsavebox{\tfive}
\sbox{\tfive}{
    \tree{\vertex1{0.6}\vertex10\vertex3{-0.6}\vertex10}
}
\caption{Coefficients of the elementary differentials with trees of the shape \usebox{\tfive} with the black root for the P-series expansions of the exact solution, the fourth order AVF collocation method, and the proposed method.}
\label{table:9-5}
{\renewcommand\arraystretch{1.5}
\begin{tabular}{cccccccccc}
$t$
& \tree{\vertex1{0.6}\vertex10\vertex3{-0.6}\vertex10}
& \tree{\white\vertex1{0.6}\black\vertex10\vertex3{-0.6}\vertex10}
& \tree{\vertex1{0.6}\white\vertex10\black\vertex3{-0.6}\vertex10}
& \tree{\white\vertex1{0.6}\vertex10\black\vertex3{-0.6}\vertex10}
& \tree{\white\vertex1{0.6}\black\vertex10\white\vertex3{-0.6}\black\vertex10}
& \tree{\white\vertex1{0.6}\black\vertex10\vertex3{-0.6}\white\vertex10}
& \tree{\vertex1{0.6}\white\vertex10\vertex3{-0.6}\vertex10}
& \tree{\white\vertex1{0.6}\black\vertex10\white\vertex3{-0.6}\vertex10}
& \tree{\white\vertex1{0.6}\vertex10\vertex3{-0.6}\vertex10} \\
exact solution &
$\tfrac{1}{20}$ &
$\tfrac{1}{20}$ &
$\tfrac{1}{20}$ &
$\tfrac{1}{20}$ &
$\tfrac{1}{20}$ &
$\tfrac{1}{20}$ &
$\tfrac{1}{20}$ &
$\tfrac{1}{20}$ &
$\tfrac{1}{20}$ \\
\hline
AVF(4) &
$\tfrac{1}{20}$ &
$\tfrac{7}{144}$ &
$\tfrac{1}{20}$ &
$\tfrac{7}{144}$ &
$\tfrac{7}{144}$ &
$\tfrac{7}{144}$ &
$\tfrac{7}{144}$ &
$\tfrac{7}{144}$ &
$\tfrac{7}{144}$ \\
\hline
proposed &
$\tfrac{1}{20}$ &
$\tfrac{\eqref{A}}{4}+\tfrac{7}{144}$ &
$\tfrac{1}{20}$ &
$\tfrac{\eqref{A}}{4}+\tfrac{7}{144}$ &
$\tfrac{\eqref{B}}{4}$ &
$\tfrac{\eqref{A}}{4}+\tfrac{7}{144}$ &
$\tfrac{\eqref{A}}{4}+\tfrac{7}{144}$ &
$\tfrac{\eqref{B}}{4}$ &
$\tfrac{\eqref{B}}{4}$ 

\end{tabular}
}
\end{table}

\begin{table}
\newsavebox{\tsix}
\sbox{\tsix}{
    \black\tree{\vertex10\vertex1{-1}\vertex2{0}\vertex3{1}}
}
\caption{Coefficients of the elementary differentials with trees of the shape \usebox{\tsix} with the black root for the P-series expansions of the exact solution, the fourth order AVF collocation method, and the proposed method.}
\label{table:9-6}
{\renewcommand\arraystretch{1.5}
\begin{tabular}{cccccccccc}
$t$ 
& \black\tree{\vertex10\vertex1{-1}\vertex2{0}\vertex3{1}}
& \black\tree{\white\vertex10\black\vertex1{-1}\vertex2{0}\vertex3{1}}
& \black\tree{\vertex10\vertex1{-1}\vertex2{0}\white\vertex3{1}}
& \black\tree{\white\vertex10\black\vertex1{-1}\vertex2{0}\white\vertex3{1}}
& \black\tree{\vertex10\vertex1{-1}\white\vertex2{0}\vertex3{1}}
& \black\tree{\vertex10\white\vertex1{-1}\vertex2{0}\vertex3{1}}
& \black\tree{\white\vertex10\vertex1{-1}\vertex2{0}\black\vertex3{1}}
& \black\tree{\white\vertex10\vertex1{-1}\vertex2{0}\vertex3{1}} \\
exact solution &
$\tfrac{1}{20}$ &
$\tfrac{1}{20}$ &
$\tfrac{1}{20}$ &
$\tfrac{1}{20}$ &
$\tfrac{1}{20}$ &
$\tfrac{1}{20}$ &
$\tfrac{1}{20}$ &
$\tfrac{1}{20}$ \\
\hline 
AVF(4) &
$\tfrac{1}{20}$ &
$\tfrac{1}{20}$ &
$\tfrac{1}{18}$ &
$\tfrac{1}{18}$ &
$\tfrac{1}{18}$ &
$\tfrac{1}{18}$ &
$\tfrac{1}{18}$ &
$\tfrac{1}{18}$ \\
\hline 
proposed &
$\tfrac{1}{20}$ &
$\tfrac{1}{20}$ &
$2\eqref{E}+\tfrac{1}{18}$ &
$2\eqref{E}+\tfrac{1}{18}$ &
$-\eqref{H} + \tfrac{1}{16}$ &
$-\eqref{B} + \tfrac{1}{4}$ &
$-\eqref{H} + \tfrac{1}{16}$ &
$-\eqref{B} + \tfrac{1}{4}$ 

\end{tabular}
}
\end{table}

\begin{table}
\newsavebox{\tseven}
\sbox{\tseven}{
    \black\tree{\vertex10\vertex1{-0.6}\vertex2{0.6}\vertex10} 
}
\caption{Coefficients of the elementary differentials with trees of the shape \usebox{\tseven} with the black root for the P-series expansions of the exact solution, the fourth order AVF collocation method, and the proposed method.}
\label{table:9-7}
{\renewcommand\arraystretch{1.5}
\begin{tabular}{ccccccccccccccccc}
$t$ 
& \black\tree{\vertex10\vertex1{-0.6}\vertex2{0.6}\vertex10} 
& \black\tree{\white\vertex10\black\vertex1{-0.6}\vertex2{0.6}\vertex10} 
& \black\tree{\vertex10\white\vertex1{-0.6}\black\vertex2{0.6}\vertex10} 
& \black\tree{\vertex10\vertex1{-0.6}\white\vertex2{0.6}\black\vertex10} 
& \black\tree{\vertex10\vertex1{-0.6}\vertex2{0.6}\white\vertex10} 
& \black\tree{\white\vertex10\vertex1{-0.6}\black\vertex2{0.6}\vertex10} 
& \black\tree{\white\vertex10\black\vertex1{-0.6}\white\vertex2{0.6}\black\vertex10} 
& \black\tree{\white\vertex10\black\vertex1{-0.6}\vertex2{0.6}\white\vertex10} 
& \black\tree{\vertex10\white\vertex1{-0.6}\vertex2{0.6}\black\vertex10} 
& \black\tree{\vertex10\white\vertex1{-0.6}\black\vertex2{0.6}\white\vertex10} 
 \\
exact solution &
$\tfrac{1}{40}$ &
$\tfrac{1}{40}$ &
$\tfrac{1}{40}$ &
$\tfrac{1}{40}$ &
$\tfrac{1}{40}$ &
$\tfrac{1}{40}$ &
$\tfrac{1}{40}$ &
$\tfrac{1}{40}$ &
$\tfrac{1}{40}$ &
$\tfrac{1}{40}$ \\
\hline
AVF(4) &
$\tfrac{1}{40}$ &
$\tfrac{1}{40}$ &
$\tfrac{1}{36}$ &
$\tfrac{1}{36}$ &
$\tfrac{1}{40}$ &
$\tfrac{1}{36}$ &
$\tfrac{1}{36}$ &
$\tfrac{1}{36}$ &
$\tfrac{1}{36}$ &
$\tfrac{1}{36}$ \\
\hline
proposed &
$\tfrac{1}{40}$ &
$\tfrac{1}{40}$ &
$\eqref{E}+\tfrac{1}{36}$ &
$-\tfrac{\eqref{H}}{2} + \tfrac{1}{32}$ &
$\tfrac{1}{40}$ &
$\eqref{E}+\tfrac{1}{36}$ &
$-\tfrac{\eqref{H}}{2} + \tfrac{1}{32}$ &
$\tfrac{1}{40}$ &
$-\tfrac{\eqref{B}}{2} + \tfrac{1}{8}$ &
$\eqref{E}+\tfrac{1}{36}$ 
\end{tabular}
}

{\renewcommand\arraystretch{1.5}
\begin{tabular}{ccccccccccccccccc}
$t$ 
& \black\tree{\vertex10\vertex1{-0.6}\white\vertex2{0.6}\vertex10} 
& \black\tree{\vertex10\white\vertex1{-0.6}\vertex2{0.6}\vertex10} 
& \black\tree{\white\vertex10\black\vertex1{-0.6}\white\vertex2{0.6}\vertex10} 
& \black\tree{\white\vertex10\vertex1{-0.6}\black\vertex2{0.6}\white\vertex10} 
& \black\tree{\white\vertex10\vertex1{-0.6}\vertex2{0.6}\black\vertex10} 
& \black\tree{\white\vertex10\vertex1{-0.6}\vertex2{0.6}\vertex10} \\
exact solution &
$\tfrac{1}{40}$ &
$\tfrac{1}{40}$ &
$\tfrac{1}{40}$ &
$\tfrac{1}{40}$ &
$\tfrac{1}{40}$ &
$\tfrac{1}{40}$ \\
\hline
AVF(4) &
$\tfrac{1}{36}$ &
$\tfrac{1}{36}$ &
$\tfrac{1}{36}$ &
$\tfrac{1}{36}$ &
$\tfrac{1}{36}$ &
$\tfrac{1}{36}$ \\
\hline
proposed &
$-\tfrac{\eqref{H}}{2} + \tfrac{1}{32}$ &
$-\tfrac{\eqref{B}}{2} + \tfrac{1}{8}$ &
$-\tfrac{\eqref{H}}{2} + \tfrac{1}{32}$ &
$-\tfrac{\eqref{B}}{2} + \tfrac{1}{8}$ &
$\eqref{E}+\tfrac{1}{36}$ &
$-\tfrac{\eqref{B}}{2} + \tfrac{1}{8}$

\end{tabular}
}
\end{table}

\begin{table}
\newsavebox{\teight}
\sbox{\teight}{
    \black\tree{\vertex10\vertex10\vertex1{-0.6}\vertex2{0.6}}
}
\caption{Coefficients of the elementary differentials with trees of the shape \usebox{\teight} with the black root for the P-series expansions of the exact solution, the fourth order AVF collocation method, and the proposed method.}
\label{table:9-8}
{\renewcommand\arraystretch{1.5}
\begin{tabular}{ccccccccccccc}
$t$ 
& \black\tree{\vertex10\vertex10\vertex1{-0.6}\vertex2{0.6}}
& \black\tree{\white\vertex10\black\vertex10\vertex1{-0.6}\vertex2{0.6}}
& \black\tree{\vertex10\white\vertex10\black\vertex1{-0.6}\vertex2{0.6}}
& \black\tree{\vertex10\vertex10\vertex1{-0.6}\white\vertex2{0.6}}
& \black\tree{\white\vertex10\vertex10\black\vertex1{-0.6}\vertex2{0.6}}
& \black\tree{\white\vertex10\black\vertex10\vertex1{-0.6}\white\vertex2{0.6}}
& \black\tree{\vertex10\white\vertex10\black\vertex1{-0.6}\white\vertex2{0.6}}
& \black\tree{\vertex10\vertex10\white\vertex1{-0.6}\vertex2{0.6}}
 \\
exact solution &
$\tfrac{1}{60}$ &
$\tfrac{1}{60}$ &
$\tfrac{1}{60}$ &
$\tfrac{1}{60}$ &
$\tfrac{1}{60}$ &
$\tfrac{1}{60}$ &
$\tfrac{1}{60}$ &
$\tfrac{1}{60}$ \\
\hline
AVF(4) &
$\tfrac{1}{72}$ &
$\tfrac{1}{72}$ &
$\tfrac{1}{72}$ &
$\tfrac{1}{72}$ &
$\tfrac{1}{72}$ &
$\tfrac{1}{72}$ &
$\tfrac{1}{72}$ &
$\tfrac{1}{72}$ \\
\hline 
proposed &
$-\eqref{C}+\tfrac{1}{12}$ &
$-\eqref{C}+\tfrac{1}{12}$ &
$-\eqref{D}+\tfrac{1}{72}$ &
$-\eqref{E}+\tfrac{1}{72}$ &
$-\eqref{D}+\tfrac{1}{72}$ &
$-\eqref{E}+\tfrac{1}{72}$ &
$\eqref{F}+\tfrac{1}{72}$ &
$\tfrac{\eqref{A}}{2}+\tfrac{1}{72}$ 

\end{tabular}
}

{\renewcommand\arraystretch{1.5}
\begin{tabular}{ccccccccccccc}
$t$ 
& \black\tree{\vertex10\white\vertex10\vertex1{-0.6}\vertex2{0.6}}
& \black\tree{\white\vertex10\black\vertex10\white\vertex1{-0.6}\vertex2{0.6}}
& \black\tree{\white\vertex10\vertex10\vertex1{-0.6}\black\vertex2{0.6}}
& \black\tree{\white\vertex10\vertex10\vertex1{-0.6}\vertex2{0.6}} \\
exact solution &
$\tfrac{1}{60}$ &
$\tfrac{1}{60}$ &
$\tfrac{1}{60}$ &
$\tfrac{1}{60}$ \\
\hline
AVF(4) &
$\tfrac{1}{72}$ &
$\tfrac{1}{72}$ &
$\tfrac{1}{72}$ &
$\tfrac{1}{72}$ \\
\hline 
proposed &
$-\eqref{G}+\tfrac{1}{72}$ &
$\tfrac{\eqref{A}}{2}+\tfrac{1}{72}$ &
$\eqref{F}+\tfrac{1}{72}$ &
$-\eqref{G}+\tfrac{1}{72}$

\end{tabular}
}
\end{table}

\begin{table}
\newsavebox{\tnine}
\sbox{\tnine}{
    \black\tree{\vertex10\vertex10\vertex10\vertex10}
}
\caption{Coefficients of the elementary differentials with trees of the shape \usebox{\tnine} with the black root for the P-series expansions of the exact solution, the fourth order AVF collocation method, and the proposed method.}
\label{table:9-9}
{\renewcommand\arraystretch{1.5}
\begin{tabular}{ccccccccccccccccc}
$t$ 
& \black\tree{\vertex10\vertex10\vertex10\vertex10}
& \black\tree{\white\vertex10\black\vertex10\vertex10\vertex10}
& \black\tree{\vertex10\white\vertex10\black\vertex10\vertex10}
& \black\tree{\vertex10\vertex10\white\vertex10\black\vertex10}
& \black\tree{\vertex10\vertex10\vertex10\white\vertex10}
& \black\tree{\white\vertex10\vertex10\black\vertex10\vertex10}
& \black\tree{\white\vertex10\black\vertex10\white\vertex10\black\vertex10}
& \black\tree{\white\vertex10\black\vertex10\vertex10\white\vertex10}
\\
exact solution &
$\tfrac{1}{120}$ &
$\tfrac{1}{120}$ &
$\tfrac{1}{120}$ &
$\tfrac{1}{120}$ &
$\tfrac{1}{120}$ &
$\tfrac{1}{120}$ &
$\tfrac{1}{120}$ &
$\tfrac{1}{120}$ \\
\hline
AVF(4) &
$\tfrac{1}{144}$ &
$\tfrac{1}{144}$ &
$\tfrac{1}{144}$ &
$\tfrac{1}{144}$ &
$\tfrac{1}{144}$ &
$\tfrac{1}{144}$ &
$\tfrac{1}{144}$ &
$\tfrac{1}{144}$ \\
\hline
proposed &
$-\tfrac{\eqref{C}}{2}+\tfrac{1}{24}$ &
$-\tfrac{\eqref{C}}{2}+\tfrac{1}{24}$ &
$-\tfrac{\eqref{D}}{2}+\tfrac{1}{144}$ &
$\tfrac{\eqref{A}}{4}+\tfrac{1}{144}$ &
$-\tfrac{\eqref{C}}{2}+\tfrac{1}{24}$ &
$-\tfrac{\eqref{D}}{2}+\tfrac{1}{144}$ &
$\tfrac{\eqref{A}}{4}+\tfrac{1}{144}$ &
$-\tfrac{\eqref{C}}{2}+\tfrac{1}{24}$ 
\end{tabular}
}

{\renewcommand\arraystretch{1.5}
\begin{tabular}{ccccccccccccccccc}
$t$ 
& \black\tree{\vertex10\white\vertex10\vertex10\black\vertex10}
& \black\tree{\vertex10\white\vertex10\black\vertex10\white\vertex10}
& \black\tree{\vertex10\vertex10\white\vertex10\vertex10}
& \black\tree{\vertex10\white\vertex10\vertex10\vertex10}
& \black\tree{\white\vertex10\black\vertex10\white\vertex10\vertex10}
& \black\tree{\white\vertex10\vertex10\black\vertex10\white\vertex10}
& \black\tree{\white\vertex10\vertex10\vertex10\black\vertex10}
& \black\tree{\white\vertex10\vertex10\vertex10\vertex10}  \\
exact solution &
$\tfrac{1}{120}$ &
$\tfrac{1}{120}$ &
$\tfrac{1}{120}$ &
$\tfrac{1}{120}$ &
$\tfrac{1}{120}$ &
$\tfrac{1}{120}$ &
$\tfrac{1}{120}$ &
$\tfrac{1}{120}$ \\
\hline
AVF(4) &
$\tfrac{1}{144}$ &
$\tfrac{1}{144}$ &
$\tfrac{1}{144}$ &
$\tfrac{1}{144}$ &
$\tfrac{1}{144}$ &
$\tfrac{1}{144}$ &
$\tfrac{1}{144}$ &
$\tfrac{1}{144}$ \\
\hline
proposed &
$-\tfrac{\eqref{G}}{2}+\tfrac{1}{144}$ &
$-\tfrac{\eqref{D}}{2}+\tfrac{1}{144}$ &
$\tfrac{\eqref{A}}{4}+\tfrac{1}{144}$ &
$-\tfrac{\eqref{G}}{2}+\tfrac{1}{144}$ &
$\tfrac{\eqref{A}}{4}+\tfrac{1}{144}$ &
$-\tfrac{\eqref{D}}{2}+\tfrac{1}{144}$ &
$-\tfrac{\eqref{G}}{2}+\tfrac{1}{144}$ &
$-\tfrac{\eqref{G}}{2}+\tfrac{1}{144}$ 
\end{tabular}
}
\end{table}

\end{document}

%% file: macro.tex
\newcommand{\bbR}{\mathbb{R}}

\newcommand{\rmd}{\mathrm{d}}

\newcommand{\trans}{{\text{\tiny\sf T}}}

\DeclareMathOperator{\diag}{diag}